\documentclass[10pt]{article}

\usepackage{eucal}

\usepackage{amsfonts,amsmath,amssymb,latexsym}

\begin{document}

\newtheorem{theorem}{Theorem}

\newtheorem{proposition}{Proposition}

\newtheorem{corollary}{Corollary}

\newtheorem{lemma}{Lemma}

\newtheorem{definition}{Definition}

\newtheorem{remark}{Remark}

\newtheorem{claim}{Claim}

\def\Ss{\mathbb{S}}

\def\Mm{\mathbb{M}^2({\kappa})}

\def\Cc{\mathcal{C}}

\def\Rr{\mathbb{R}}

\def\Pp{\mathbb{P}_t}

\def\dd{\mathrm{d}}

\def\mt{\mathcal}

\def\ae{\langle}

\def\ad{\rangle}

\def\sn{\textrm{sn}}

\def\ct{\textrm{ct}}

\def\cs{\textrm{cs}}

\def\re{\textrm{Re}}

\def\im{\textrm{Im}}

\title{Examples and structure of CMC surfaces in some Riemannian
and Lorentzian  homogeneous spaces}

\author{Marcos P. de A. Cavalcante\thanks{Partially supported by CNPq.} \,
and \, Jorge H. S. de Lira\thanks {Partially supported by CNPq   and
FUNCAP.}}

\date{}

\maketitle

\begin{abstract}
It is proved that the holomorphic quadratic differential associated
to CMC surfaces in Riemannian products $\Ss^2\times\Rr$ and
$\mathbb{H}^2\times \Rr$ discovered by U. Abresch and H. Rosenberg
could be obtained as a linear combination of usual Hopf
differentials. Using this fact, we are able to extend it for
Lorentzian products. Families of examples of helicoidal CMC surfaces
on these spaces are explicitly described. We also present some
characterizations of CMC rotationally invariant discs and spheres.
Finally, after establish some height and area estimates, we prove
the existence of constant mean curvature Killing graphs.
\end{abstract}
\vspace{0.3cm}

{\small \noindent {\bf Keywords:} constant mean curvature,
holomorphic quadratic differentials, Killing graphs

\noindent {\bf MSC 2000:} 53C42, 53A10.}

\section{Introduction}

U. Abresch and H. Rosenberg had recently proved that there exists a
quadratic differential for an  immersed surface in
$\mathbb{M}^2(\kappa)\times \Rr$ which is holomorphic when the
surface has constant mean curvature. Here, $\mathbb{M}^2(\kappa)$
denotes the two-dimensional simply connected space form with
constant curvature $\kappa$.  This differential $Q$ plays the role
of the usual Hopf differential in the theory of constant mean
curvature surfaces immersed in space forms. Thus, they were able to
prove the following theorem: \vspace{0.3cm}

\noindent{\bf Theorem.} (Theorem 2, p. 143, \cite{AR}) {\it Any
immersed cmc sphere $S^2\looparrowright \mathbb{M}^2(\kappa)\times
\Rr$ in a product space is actually one of the embedded rotationally
invariant cmc spheres $S^2_H\subset \mathbb{M}^2(\kappa) \times
\Rr$.} \vspace{0.3cm}

The rotationally invariant spheres referred to above were
constructed independently by W.-Y. Hsiang and W.-T. Hsiang in
\cite{HH} and by R. Pedrosa and M. Ritor\'e in \cite{P} and
\cite{PR}. The theorem quoted above proves affirmatively a
conjecture stated by Hsiang and Hsiang in their paper \cite{HH}.
More importantly, it indicates that some tools often used for
surface theory in space forms could be redesigned to more general
three dimensional homogeneous spaces, the more natural ones after
space forms being $\mathbb{M}^2(k)\times \Rr$. The price to be paid
in abandoning space forms is that the technical difficulties are
more involved. The method in \cite{AR} is to study very closely the
revolution surfaces in $\mathbb{M}^2(\kappa)\times \Rr$ in order to
guess the suitable differential.

Our idea here is to relate the $Q$ differential on a surface
$\Sigma$ immersed in $\Mm\times\Rr$ with the usual Hopf differential
after embedding $\mathbb{M}^2(\kappa) \times \Rr$ in some Euclidean
space $\mathbb{E}^4$. We prove that $Q$ is written as a linear
combination of the Hopf differentials $\Psi^1$ and $\Psi^2$
associated to two normal directions spanning the normal bundle of
$\Sigma$ in $\mathbb{E}^4$. This fact is also true when the product
$\Mm\times \Rr$ carries a Lorentzian metric. More precisely, if we
define $r$ as $r^2 =\epsilon/\kappa$ for $\epsilon=\textrm{sgn}\,
\kappa$ we state the following result:

\vspace{0.3cm} \noindent{\bf Theorem. }{(Theorem 7, p. 25)  \it The
quadratic differential $Q=2H \Psi^1-\varepsilon\frac{\epsilon}{r}\,
\Psi^2$ is holomorphic on $\Sigma\looparrowright\Mm\times\Rr$ if the
mean curvature $H$ of $\Sigma$ is constant. Inversely, if we suppose
that $\Sigma$ is compact (more generally, if $\Sigma$ does not admit
a function without critical points, or a vector field without
singularities), then $H$ is constant if $Q$ is holomorphic.}

\vspace{0.3cm}

Our aim here is to explore geometrical consequences of this
alternative presentation of $Q$. We next give a brief description of
this paper. The sections 2 and 3 are concerned with the existence
and structure of families of isometric surfaces with same constant
mean curvature on both Riemannian and Lorentzian products which are
invariant by certain isometry groups of the ambient space. Our
construction is inspired by that one presented in \cite{dCD} and
\cite{ET}. In Section 4, we present the proof of the Theorem 7 and a
variant of the classical Theorem of Joachimstahl which gives a
characterization of CMC rotationally invariant  discs and spheres in
the same spirit of the result by Abresch and Rosenberg mentioned
above (see Theorem 8).

We also prove on Section 5  the following result about free boundary
CMC surfaces, based on  the well-known Nitsche's work on
partitioning problem:

\vspace{0.3cm}

\noindent {\bf Theorem. }{(Theorem 9, p. 29) \it Let $\Sigma$ be a
surface immersed in $\Mm\times \Rr$ whose boundary is contained in
some horizontal plane $\mathbb{P}_a$. Suppose that $\Sigma$ has
constant mean curvature and that its angle with $\mathbb{P}_a$ is
constant along its boundary.  If $\varepsilon =1$ and $\Sigma$ is
disc-type, then $\Sigma$ is a spherical cap. If $\varepsilon =-1$,
then $\Sigma$ is a hyperbolic cap.}

\vspace{0.3cm}

The variational meaning of the conditions on $\Sigma$  could be seen
on Section 5. We end this section with a characterization of stable
CMC discs with circular boundary on $\Mm\times \Rr$ which
generalizes a nice result of Al\'ias, L\'opez and Palmer (see
\cite{ALP}). Finally, on Section 6, we obtain estimates of some
geometrical data of CMC surfaces with boundary lying on vertical
planes in $\Mm\times \Rr$. These estimates are then used to prove
the existence on non-negatively curved Riemannian products of CMC
Killing graphs with boundary contained in vertical planes:

\vspace{0.3cm} \noindent {\bf Theorem. }{(Theorem 12, p. 35) \it Let
$\Pi$ be a vertical plane on the Riemannian product $\Mm\times\Rr$,
$\kappa\ge 0$, determined by an unit vector $a$ in $\mathbb{E}^4$.
Let $\Omega$ be a domain on $\Pi$ which does not contain points of
the axis $\{\pm a\}\times\Rr$. If $|H|< \kappa_g/\tilde\gamma$,
where $\kappa_g$ is the geodesic curvature of $\partial\Omega$ in
$\Pi$, then there exists a surface (a Killing graph) with constant
mean curvature $H$ and boundary $\partial\Omega$.}

\vspace{0.3cm}

The constant $\tilde\gamma$ depends on the maximum and minimum
values on $\Omega$ of the norm of the Killing vector field generated
by rotations fixing $a$.

In a forthcoming paper (see \cite{L}), one of the authors elaborates
versions of the results contained here for constant mean curvature
hypersurfaces in  some homogeneous  spaces and warped products.
There, a suitable treatment of Minkowski formulae  gives some hints
about stability problems and the existence of general Killing
graphs.

\vspace{.3cm}

\noindent {\bf Acknowledgments:} The first author acknowledges the 
hospitality of the Departamento de Matem\'atica of Universidade Federal 
do Cear\'a in the Summer of 2005.

\section{Screw-motion invariant CMC surfaces}

\subsection{The mean curvature equation}

Let $\Mm$ be a two dimensional simply connected surface endowed with
a Riemannian  complete metric $\dd\sigma^2$ with constant sectional
curvature $\kappa$. We fix the metric $\varepsilon\dd t^2 +
\dd\sigma^2$, $\varepsilon=\pm 1$, on the product $\Mm\times \Rr$.
This metric is Lorentzian if  $\varepsilon = -1$ and Riemannian if
$\varepsilon = 1$.

A tangent vector $v$ to $\Mm\times \Rr$ is projected on horizontal
component $v^h$ and vertical component $v^t$, respectively tangent
to the $T\Mm$ and $T\Rr$ factors. We denote  by
$\langle\cdot,\cdot\rangle$ and $D$ respectively  the metric and
covariant derivative in $\Mm\times \Rr$. The curvature tensor  in
$\Mm\times\Rr$ is denoted by $\bar R$.

Let $(\rho, \theta)$ be polar coordinates centered at some point
$p_0$ in $\Mm$ and the corresponding cylindrical coordinates $(\rho,
\theta, t)$ in $\Mm\times \Rr$. Fix then a curve $s\mapsto
(\rho(s),0, t(s))$ in the plane  $\theta=0$. If we rotate this curve
at the same time we translate it along the $t$ axis with constant
speed $b$, we obtain a screw-motion invariant surface (for short, an
{\it helicoidal surface}) $\Sigma$ in $\Mm\times \Rr$ whose axis is
$\{p_0\}\times \Rr$. This means that this surface has a
parametrization $X$, in terms of the cylindrical cordinates defined
above, of the following form:
\begin{equation}
\label{X} X(s, \theta)= (\rho(s), \theta, t(s) + b\, \theta).
\end{equation}

For $b=0$ the surface $\Sigma$ is a {\it revolution surface}, i.e.,
it is invariant with respect to the action of $O(2)$ on $\Mm\times
\Rr$ fixing the axis $\{p_0\}\times \Rr$. Another interesting
particular case is obtained when $t(s)= 0$ and $s\mapsto \rho (s)$
is just  an arbitrary parametrization of the horizontal geodesic
$\theta=0, \, t=0$. Here, the resulting surfaces are called {\it
helicoids}. We will see that helicoids are examples with zero mean
curvature. Helicoidal surfaces into Riemannian products $\Mm\times
\Rr$ were already extensively studied in \cite{dCD}, \cite{PR},
\cite{HH}, \cite{AR}, \cite{ET} and \cite{MO}, for instance. In
Lorentzian products, we will consider only {\it space-like}
helicoidal surfaces, i.e., surfaces for which the metric induced on
them is a Riemannian metric.

The tangent plane to $\Sigma$ at a point $X(s, \theta)$ is spanned
by the coordinate vector fields
\begin{eqnarray*}
X_s = \dot \rho \,\partial_\rho + \dot t \, \partial_t \quad
\textrm{and} \quad X_\theta = \partial_\theta + b\, \partial_t.
\end{eqnarray*}
Throughout this text, we denote $\sn_\kappa(\rho)=|\ae
\partial_\theta, \partial_\theta\ad|^{1/2}$.  For further reference, we still
denote $\cs_\kappa(\rho)=:\frac{\dd}{\dd\rho}\sn_\kappa(\rho)$. With
this notation, an orientation for $X(\Sigma)$ is given by the unit
normal vector field
\[
n = \frac{1}{W}\, \big( \sn_\kappa(\rho)\, \dot t \,\partial_\rho +
b\, \sn^{-1}_\kappa(\rho)\,\dot \rho \,\partial_\theta -\varepsilon
\sn_\kappa(\rho)\,\dot \rho\,
\partial_t\big),
\]
where
\[
W^2 = :\sn^2_\kappa (\rho) (\dot \rho^2+\varepsilon \dot t^2)
+\varepsilon b^2 \dot \rho^2.
\]
We suppose that $\langle n,n\rangle =\varepsilon$. When
$\varepsilon=-1$ this assumption implies that $\Sigma$ is space-like
and  that $W^2>0$. It also follows  that $U^2=:
\sn_\kappa^2(\rho)+\varepsilon b^2>0$. The induced metric on
$\Sigma$ is given by
\begin{eqnarray*}
\ae\dd X,  \dd X \ad & = & E \dd s^2 +2F\dd s \dd \theta + G\dd
\theta^2\\
& =: &  \big(\dot \rho^2 +\varepsilon\dot t^2\big)\, \dd s^2
+2\varepsilon b\dot t \, \dd s \,\dd \theta  +
\big(\sn^2_\kappa(\rho) +\varepsilon b^2\big)\, \dd \theta^2.
\end{eqnarray*}
The vector field $\partial_t$ is parallel and $s\mapsto \rho(s)$
parametrizes a geodesic on $\Mm$. So
 it follows that
\[
X_{ss}= :D_{X_s}X_s = \ddot \rho \partial_\rho + \ddot t \partial_t.
\]
The remaining two covariant derivatives of the coordinate vector
fields are
\[
X_{s\theta}=:D_{X_s}X_\theta = \dot \rho\,
D_{\partial_\rho}\partial_\theta, \quad X_{\theta\theta} =:
D_{X_\theta}X_\theta = D_{\partial_\theta}\partial_\theta.
\]
The first coefficient of the second fundamental form $-\ae \dd n,
\dd X\ad$ of $\Sigma$ is given by
\[
e=:\ae X_{ss},n\ad = \frac{\sn_\kappa(\rho)}{W}\,(\ddot \rho \dot t-
\ddot t \dot\rho)
\]
and since  $\ae D_{\partial_\rho}\partial_\theta,
\partial_\rho\ad = \ae
\partial_\theta, D_{\partial_\rho}\partial_\rho\ad=0$ and $\ae
D_{\partial_\rho}\partial_\theta, \partial_\theta\ad =
\frac{1}{2}\frac{\dd}{\dd \rho}
\sn_\kappa^2(\rho)=\sn_\kappa(\rho)\cs_\kappa(\rho)$ it follows that
\begin{eqnarray*}
& & f=: \ae  X_{s\theta}, n \ad  = \frac{1}{W}\, b \dot \rho^2
\cs_\kappa (\rho).
\end{eqnarray*}
Finally,   $\ae D_{\partial_\theta}\partial_\theta,
\partial_\theta \ad =0$ implies
\begin{eqnarray*} & & g=:\ae X_{\theta\theta} , n\ad =
 -\frac{1}{W}\sn_\kappa(\rho) \dot t \, \ae
\partial_\theta, D_{\partial_\rho}\partial_\theta\ad=
-\frac{1}{W}\dot t\,\sn_\kappa^2(\rho) \cs_\kappa(\rho).
\end{eqnarray*}
Thus the formula $H=\frac{1}{2}\, (eG-2fF+gE)/(EG-F^2)$ for the mean
curvature of $X$ reads
\begin{eqnarray}
\label{H-eq} & & 2HW^3 =(\ddot \rho \dot t-\ddot t \dot
\rho)\sn_\kappa(\rho)(\sn_\kappa^2(\rho)+\varepsilon
b^2)-2\varepsilon b^2\dot t \dot \rho^2
\cs_\kappa(\rho)\nonumber\\
& &-\dot t (\dot \rho^2+\varepsilon \dot
t^2)\sn_\kappa^2(\rho)\cs_\kappa(\rho).
\end{eqnarray}
We suppose  momentarily that the profile curve $(\rho(s),0,t(s))$ is
given as a graph $t=t(\rho)$. Thus we put $\rho=s$ above and find
\begin{equation}
\label{W} W^2= EG-F^2= \sn_\kappa^2(\rho)(1+\varepsilon\dot
t^2)+\varepsilon b^2.
\end{equation}
Therefore the mean curvature equation (\ref{H-eq}) reduces to
\begin{equation}
\label{H-eq2} 2HW^3= -\ddot t
\sn_\kappa(\rho)(\sn_\kappa^2(\rho)+\varepsilon b^2)-2\varepsilon
b^2\dot t \dot \sn_\kappa(\rho)-\dot t (1+\varepsilon\dot
t^2)\sn_\kappa^2(\rho)\dot \sn_\kappa(\rho),
\end{equation}
where the derivatives are taken with respect to the parameter
$\rho$. One easily verifies that the expression
\[
\frac{\dd}{\dd \rho}\Big(\frac{\dot
t\sn_\kappa^2(\rho)}{W}\Big)=-2H\sn_\kappa(\rho)
\]
is equivalent to the equation (\ref{H-eq2}) above. This means that
\begin{equation}
\label{integral} \frac{\dd t}{\dd \rho}\frac{\sn_\kappa^2(\rho)}{W}
= I-2H\int\sn_\kappa(\rho)\, \dd \rho
\end{equation}
is a first integral to the mean curvature equation (\ref{H-eq2})
associated to translations on $t$ axis.

\subsection{A Bour's type lemma and rotational examples}

Next, we will obtain orthogonal parameters for $\Sigma$ for which
one  of the families of coordinate curves is given by geodesics on
$\Sigma$. For this, we write
\begin{eqnarray*}
&  & \ae \dd X, \dd X\ad=   \big(\dot \rho^2 +\varepsilon\dot
t^2\big)\, \dd s^2 + (\sn_\kappa^2(\rho) +\varepsilon b^2)\big(\dd
\theta + U^{-2}\varepsilon b\dot t \dd s\big)^2 - U^{-2} b^2 \dot
t^2\, \dd s^2\\
& &= \frac{W^2}{U^2}\, \dd s^2 + U^2 \dd \tilde\theta^2 = \dd \tilde
s^2 + U^2 \dd \tilde\theta^2,
\end{eqnarray*}
where $\dd \tilde s = \frac{W}{U}\, \dd s$ and $\dd\tilde\theta =\dd
\theta + U^{-2}\varepsilon b \dot t \dd s$. These differentials
could be locally integrated and  furnish an actual change of
coordinates on $\Sigma$. For revolution surfaces (i.e., for $b=0$)
such change of variables is not necessary. More precisely, it
consists only in to assume that $s$ is the arc lenght of the profile
curve $s \mapsto (\rho(s),0, t(s))$. For helicoids we have $\dot t
=0$ and then the change of variables is again useless since here we
may choose $\rho =s$ along the rules of the helicoid. Since that $W$
and $U$ depend only on $s$, then $\tilde s$ is a function of $s$
only with $\frac{\dd \tilde s}{\dd s}= \frac{W}{U}$. Notice that
\[
\frac{W^2}{U^2} =\frac{\sn_\kappa^2(\rho)(\dot
\rho^2+\varepsilon\dot t^2)+\varepsilon
b^2\dot\rho^2}{\sn_\kappa^2(\rho)+\varepsilon b^2}=\dot\rho^2
+\frac{\varepsilon\sn_\kappa^2(\rho)\dot
t^2}{\sn_\kappa^2(\rho)+\varepsilon b^2}.
\]
Thus  the functions $\tilde s, \, \tilde\theta$ satisfy the system
\begin{eqnarray}
\label{eq1-system} & & \dd \tilde s^2 =\dd\rho^2
+\frac{\varepsilon\sn_\kappa^2(\rho)}{\sn_\kappa^2(\rho)+\varepsilon
b^2}\,\dd
t^2,\\
\label{eq2-system} & & U\,\dd \tilde \theta =
(\sn_\kappa^2(\rho)+\varepsilon b^2)^{1/2}\big(\dd \theta
+\frac{\varepsilon b}{\sn_\kappa^2(\rho)+\varepsilon b^2}\, \dd
t\big).
\end{eqnarray}
One easily verifies that the coordinate curves
$\tilde\theta=\textrm{cte.}$ are geodesics on $\Sigma$. In fact, if
we consider the frame $e_1=\partial_{\tilde s}$ and
$e_2=U^{-1}\partial_{\tilde \theta}$ and the associated co-frame
$\omega^1=\dd \tilde s$ and $\omega^2=U \dd \tilde \theta$, then
$\omega^2_1 =\frac{\dot U}{U}\, \omega^2$. So, if $\nabla$ denotes
the induced connection on $\Sigma$ then $\nabla_{e_1}
e_1=\nabla_{\partial_{\tilde s}}{\partial_{\tilde s}}=0$. These
geodesics intersect orthogonally the curves $\tilde
s=\textrm{cte.}$. This allows us also to prove that the intrinsic
Gaussian curvature $K_{\textrm{int}}$ of $\Sigma$ is simply
$-\frac{\ddot U}{U}$.

Now, given the (natural) parameters $(\tilde s, \tilde \theta)$ on
$\Sigma$  and the function $U(\tilde s)$ we want to determine a
two-parameter family of isometric immersions $X_{m,b}:\Sigma\to
\Mm\times \Rr$ in such a way that the immersed surfaces
$X_{m,b}(\Sigma)$ are helicoidal and have induced metric given by
$\dd \tilde s^2+ U^2 \dd \tilde \theta^2$. Moreover, we require that
the original immersion $X$ belongs to that family. For this, it
suffices that the equations (\ref{eq1-system}) and
(\ref{eq2-system}) are satisfied by coordinates $\rho, \theta, t$ as
functions of $\tilde s,\tilde\theta$ for some positive constant $b$.
We refer in what follows to the original immersion and its pitch by
$X_0$ and $b_0$.

From equations (\ref{eq1-system}) and (\ref{eq2-system}) we have
$\frac{\partial\rho}{\partial \tilde \theta}=\frac{\partial
t}{\partial \tilde \theta}=0$ and
\begin{eqnarray}
\label{partialtheta} \frac{\partial \theta}{\partial \tilde
s}=-\frac{\varepsilon b}{\sn_\kappa^2(\rho)+\varepsilon b^2}\,
\frac{\dd t}{\dd \tilde s}, \quad \frac{\partial \theta}{\partial
\tilde \theta} = \frac{U}{(\sn_\kappa^2(\rho)+\varepsilon
b^2)^{1/2}}
\end{eqnarray}
and therefore
\[
\frac{\partial^2 \theta}{\partial \tilde s\partial\tilde \theta}
=\frac{\partial^2 \theta}{\partial \tilde \theta\partial\tilde s} =
0.
\]
Hence $\frac{\partial \theta}{\partial \tilde
\theta}=U(\sn_\kappa^2(\rho)+\varepsilon b^2)^{-1/2}$ does not
depend on $\tilde s$. Since  $U(\sn_\kappa^2(\rho)+\varepsilon
b^2)^{-1/2}$ does not depend also on $\tilde \theta$ it follows that
\begin{eqnarray}
\label{U} \frac{U}{(\sn_\kappa^2(\rho)+\varepsilon
b^2)^{1/2}}=\frac{1}{m}
\end{eqnarray}
for some non zero constant $m$. This defines the first parameter of
the family. The other one is the varying pitch $b$. We have
$X_0=X_{1,b_0}$. Differentiating $m^2 U^2 = \sn_\kappa^2(\rho)
+\varepsilon b^2$ with respect to $\tilde s$ we find
\[
m^2 U \dot U = \sn_\kappa(\rho) \cs_\kappa(\rho)\dot \rho.
\]
Thus since  $\sn_\kappa^2(\rho) = m^2U^2 -\varepsilon b^2$ and
$\cs_\kappa(\rho)^2 +\kappa\sn_\kappa^2(\rho)=1$ it is clear that
\[
\cs_\kappa^2(\rho)=1-\kappa(m^2U^2-\varepsilon b^2).
\]
The differential equation for $\rho$ is then
\begin{eqnarray}
\label{dotrho} \dot \rho^2 = \frac{m^4 U^2 \dot
U^2}{(m^2U^2-\varepsilon b^2)(1-\kappa(m^2U^2-\varepsilon b^2))}.
\end{eqnarray}
From the equation (\ref{eq1-system})  we conclude that $t$ satisfies
the equation
\begin{eqnarray}
\label{dott} \dot t^2 = \frac{ m^2 U^2}{(m^2U^2-\varepsilon
b^2)^2}\Big(\varepsilon\,\frac{(m^2U^2-\varepsilon
b^2)(1-\kappa(m^2U^2-\varepsilon b^2))-m^4 U^2 \dot
U^2}{1-\kappa(m^2U^2-\varepsilon b^2)}\Big).
\end{eqnarray}
Finally we infer from (\ref{partialtheta}) and (\ref{U}) that
\begin{eqnarray}
\label{dottheta} \big(\frac{\partial \theta}{\partial \tilde
s}\big)^2 = \varepsilon b^2  \,\frac{(m^2U^2-\varepsilon
b^2)(1-\kappa(m^2U^2-\varepsilon b^2))-m^4 U^2 \dot
U^2}{m^2U^2(m^2U^2-\varepsilon b^2)^2(1-\kappa(m^2U^2-\varepsilon
b^2))}
\end{eqnarray}
and $\frac{\partial \theta}{\partial\tilde \theta}=\frac{1}{m}$.
Integrating these equations we obtain
\begin{eqnarray}
\label{coordinaterho} \rho(\tilde s)  &=& \int\bigg(\frac{m^4 U^2
\dot U^2}{(m^2U^2-\varepsilon b^2)(1-\kappa(m^2U^2-\varepsilon
b^2))}\bigg)^{1/2}\, \dd \tilde
s,\\
\label{coordinatet} t(\tilde s) &=& \int\bigg(\varepsilon
\,\frac{(m^2U^2-\varepsilon b^2)(1-\kappa(m^2U^2-\varepsilon
b^2))-m^4 U^2 \dot U^2}{1-\kappa(m^2U^2-\varepsilon
b^2)}\bigg)^{1/2}\cdot\\ \nonumber& & \frac{ m U}{m^2U^2-\varepsilon
b^2}\, \dd
\tilde s,\\
\label{coordinatetheta} \theta(\tilde s, \tilde \theta) &=&
\frac{1}{m}\, \tilde
\theta + \int \frac{ b}{mU(m^2U^2-\varepsilon b^2)}\cdot\\
\nonumber& & \bigg(\varepsilon\,\frac{(m^2U^2-\varepsilon
b^2)(1-\kappa(m^2U^2-\varepsilon b^2))-m^4 U^2 \dot
U^2}{1-\kappa(m^2U^2-\varepsilon b^2)}\bigg)^{1/2}\, \dd \tilde s,
\end{eqnarray}
with $\sn_\kappa^2(\rho) = m^2 U^2 -\varepsilon b^2$.

\vspace{0.3cm}

\noindent {\bf Theorem 1.} {\it Given a helicoidal surface $X_0:
\Sigma \to\Mm\times \Rr$, with pitch $b_0$, there exists a
two-parameter family of isometric helicoidal surfaces parametrized
by $X_{m,b}:\Sigma\to \Mm\times \Rr$ with pitch $b$ such that
$X_0=X_{1,b_0}$ with coordinates given by
(\ref{coordinaterho})-(\ref{coordinatetheta}).}

\vspace{0.3cm}

We now calculate the components of the second fundamental form and
the mean curvature of these surfaces with respect to the parameters
$(\tilde s, \tilde \theta)$. Under the change of parameters
$(s,\theta)\mapsto (\tilde s, \tilde \theta)$ the second fundamental
form becomes
\begin{eqnarray}
\label{2ndff} \nonumber& & -\ae \dd n, \dd X\ad = \Big(e
\big(\frac{\partial s}{\partial \tilde s}\big)^2+ 2f\frac{\partial
s}{\partial \tilde s}\frac{\partial \theta}{\partial \tilde s} + g
\big(\frac{\partial \theta}{\partial \tilde s}\big)^2\Big)\, \dd
\tilde s^2+ 2\,\Big(e \frac{\partial s}{\partial \tilde
s}\frac{\partial s}{\partial \tilde \theta}+ \\\nonumber & &
f\big(\frac{\partial s}{\partial \tilde s}\frac{\partial
\theta}{\partial \tilde \theta}+\frac{\partial s}{\partial \tilde
\theta}\frac{\partial \theta}{\partial \tilde s}\big) + g
\frac{\partial \theta}{\partial \tilde s}\frac{\partial
\theta}{\partial \tilde \theta}\Big)\, \dd \tilde s \dd \tilde
\theta +   \Big(e \big(\frac{\partial s}{\partial \tilde
\theta}\big)^2+ 2f\frac{\partial s}{\partial \tilde
\theta}\frac{\partial \theta}{\partial \tilde \theta} + g
\big(\frac{\partial \theta}{\partial \tilde \theta}\big)^2\Big)\,
\dd \tilde \theta^2\\
& & =: \tilde e \, \dd \tilde s^2 + 2 \tilde f \, \dd \tilde s
\tilde \theta + \tilde g \, \dd \tilde \theta^2.
\end{eqnarray}
If we choose $s=\rho$, then we have from the expressions
(\ref{coordinaterho})-(\ref{coordinatetheta}) above that
\[
\frac{\partial s}{\partial \tilde s}= \frac{\dd \rho}{\dd \tilde s},
\, \frac{\partial s}{\partial \tilde \theta}=0, \, \frac{\partial
\theta}{\partial \tilde \theta}=\frac{1}{m}.
\]
Turning back to the expression (\ref{2ndff}) one finds
\begin{eqnarray*}
& & \tilde f  =f\frac{\dd \rho}{\dd \tilde s}\frac{1}{m} + g
\frac{\partial \theta}{\partial \tilde
s}\frac{1}{m}=\frac{1}{m}\frac{1}{W}b \cs_\kappa(\rho)\frac{\dd
\rho}{\dd \tilde s}-\frac{1}{m}\frac{1}{W}\sn_\kappa^2(\rho)
\frac{\dd t}{\dd \rho}\cs_\kappa(\rho)\frac{\partial
\theta}{\partial \tilde s}
\end{eqnarray*}
and
\begin{eqnarray*}
& & \tilde g = g \big(\frac{\partial \theta}{\partial \tilde
\theta}\big)^2 =
g\frac{1}{m^2}=-\frac{1}{m^2}\frac{1}{W}\sn_\kappa^2(\rho) \frac{\dd
t}{\dd \rho}\cs_\kappa(\rho).
\end{eqnarray*}
However it holds that
\[
\frac{1}{W}\frac{\dd t}{\dd \rho}= \frac{1}{mU}\, \frac{\dd t}{\dd
\tilde s}.
\]
Thus the expressions
\[
\frac{\dd t}{\dd \tilde s}=\frac{mU}{m^2U^2-\varepsilon
b^2}\,\Big(\varepsilon\,\frac{(m^2U^2-\varepsilon
b^2)(1-\kappa(m^2U^2-\varepsilon b^2))- m^4 U^2 \dot
U^2}{1-\kappa(m^2U^2-\varepsilon b^2)}\Big)^{1/2}
\]
and $\sn_\kappa^2(\rho)=m^2U^2-\varepsilon b^2$ imply that
\begin{eqnarray}
\label{lhsI} \sn_\kappa^2(\rho)\frac{1}{W}\frac{\dd t}{\dd
\rho}=\Big(\varepsilon\,\frac{(m^2U^2-\varepsilon
b^2)(1-\kappa(m^2U^2-\varepsilon b^2))- m^4 U^2 \dot
U^2}{1-\kappa(m^2U^2-\varepsilon b^2)}\Big)^{1/2}.
\end{eqnarray}
Notice that this expression is the left-hand side of the first
integral (\ref{integral}). Thus we obtain
\begin{eqnarray}
\label{integral-final} \sqrt{\varepsilon\,\frac{(m^2U^2-\varepsilon
b^2)(1-\kappa(m^2U^2-\varepsilon b^2))- m^4 U^2 \dot
U^2}{1-\kappa(m^2U^2-\varepsilon b^2)}} =I-2H\int \sn_\kappa(\rho)
\dd \rho.
\end{eqnarray}
Since  $\cs_\kappa(\rho)=(1-\kappa(m^2U^2-\varepsilon b^2))^{1/2}$
then
\[
\tilde g = -\frac{1}{m^2}\,\sqrt{\varepsilon\big((m^2U^2-\varepsilon
b^2)(1-\kappa(m^2U^2-\varepsilon b^2))- m^4 U^2 \dot U^2\big)}.
\]
Now we calculate $\tilde f$ using (\ref{partialtheta})
\begin{eqnarray*}
& & \tilde f =\frac{b}{m^4U^3}\cs_\kappa(\rho)
\Big(\sn_\kappa^2(\rho) \big(\frac{\dd \rho}{\dd \tilde s}\big)^2
+\varepsilon b^2 \big(\frac{\dd \rho}{\dd \tilde s}\big)^2 +
\varepsilon \sn_\kappa^2(\rho)\big(\frac{\dd t}{\dd \tilde
s}\big)^2\Big)  \\
& & = \frac{b}{m^4U^3}\cs_\kappa(\rho)m^2 U^2=
\frac{b}{m^2U}\sqrt{1-\kappa(m^2U^2-\varepsilon b^2)}.
\end{eqnarray*}
Finally we calculate $\tilde e$. For this one uses the Gauss formula
$K_{\textrm{int}} - \bar K = K_{\textrm{ext}}$. Here $\bar K$ is the
ambient sectional curvature and, by definition,
$K_{\textrm{ext}}=(\tilde e\tilde f - \tilde g^2)/ U^2$. So
\begin{eqnarray*}
& & \bar K = \frac{\ae \bar R(\partial_{\tilde s},
\partial_{\tilde \theta})\partial_{\tilde s}, \partial_{\tilde \theta}\ad}
{|\partial_{\tilde s}|^2|\partial_{\tilde\theta}|^2-\ae
\partial_{\tilde s}, \partial_{\tilde
\theta}\ad^2}=\frac{\kappa}{U^2}\, \big(U^2-\varepsilon U^2\ae
\partial_{\tilde s},\partial_t \ad^2 - \varepsilon\ae
\partial_{\tilde \theta},\partial_t \ad^2\big)
\end{eqnarray*}
However equations (\ref{X}) and (\ref{partialtheta}) show that
\begin{eqnarray}
\label{tildes} \nonumber& & \ae\partial_{\tilde s},\partial_t
\ad=\varepsilon \big(\frac{\dd t}{\dd \tilde s}+ b \frac{\partial
\theta}{\partial \tilde s}\big)=\frac{\varepsilon}{m^2U^2} \big(m^2
U^2 -\varepsilon b^2\big)\cdot\frac{\dd t}{\dd \tilde s}.
\end{eqnarray}
One also finds
\begin{eqnarray}
\label{tildetheta} & & \ae \partial_{\tilde\theta}, \partial_t \ad =
\varepsilon \, \frac{\partial\theta}{\partial \tilde \theta} =
\varepsilon b \frac{1}{m}.
\end{eqnarray}
Then
\begin{eqnarray*}
& & U^2-\varepsilon U^2\ae
\partial_{\tilde s},\partial_t \ad^2 - \varepsilon\ae
\partial_{\tilde \theta},\partial_t \ad^2= U^2 -\frac{\varepsilon}{m^4U^2}(m^2U^2-\varepsilon
b^2)^2\, \big(\frac{\dd t}{\dd \tilde s}\big)^2-\frac{\varepsilon
b^2}{m^2}.
\end{eqnarray*}
Finally $K_{\textrm{int}}=-\frac{\ddot U}{U}$ yields
\begin{eqnarray*}
\tilde e \tilde g - \tilde f^2 =-U\ddot U -\kappa\big(U^2
-\frac{\varepsilon}{m^4U^2}(m^2U^2-\varepsilon b^2)^2\,
\big(\frac{\dd t}{\dd \tilde s}\big)^2-\frac{\varepsilon
b^2}{m^2}\big).
\end{eqnarray*}
Thus
\begin{eqnarray*}
& & \tilde e \tilde g =-U\ddot U -\kappa U^2 +\kappa\frac{1}{m^2}
\Big(m^2U^2-\varepsilon b^2-\frac{m^4 U^2 \dot
U^2}{1-\kappa(m^2U^2-\varepsilon
b^2)}\Big)\\
& & +\kappa\frac{\varepsilon b^2}{m^2} +\frac{b^2}{m^4
U^2}\big(1-\kappa(m^2U^2-\varepsilon b^2)\big)\\
& & =-U\ddot U  -\frac{\kappa m^2 U^2}{1-\kappa(m^2U^2-\varepsilon
b^2)}\, \dot U^2+\frac{b^2}{m^4 U^2}\big(1-\kappa(m^2U^2-\varepsilon
b^2)\big).
\end{eqnarray*}
So
\begin{eqnarray*}
\tilde e = \frac{m^2U \ddot U+\frac{\kappa m^4
U^2}{1-\kappa(m^2U^2-\varepsilon b^2)}\, \dot U^2-\frac{b^2}{m^2
U^2}\big(1-\kappa(m^2U^2-\varepsilon
b^2)\big)}{\sqrt{\varepsilon\big((m^2U^2-\varepsilon
b^2)(1-\kappa(m^2U^2-\varepsilon b^2))- m^4 U^2 \dot U^2\big)}}.
\end{eqnarray*}
The mean curvature is expressed in parameters $(\tilde s, \tilde
\theta)$ as $2H = \tilde e + \frac{\tilde g}{U^2}$. Thus we have
\begin{eqnarray*}
& & 2H \,R= m^2U \ddot U+\frac{\kappa m^4
U^2}{1-\kappa(m^2U^2-\varepsilon b^2)}\, \dot U^2-\frac{b^2}{m^2
U^2}\big(1-\kappa(m^2U^2-\varepsilon b^2)\big) \\
& &-\frac{1}{m^2U^2}\, R^2 = m^2U \ddot U+\big(m^2 +\frac{\kappa m^4
U^2}{1-\kappa(m^2U^2-\varepsilon b^2)}\big)\, \dot U^2- \varepsilon
\big(1-\kappa (m^2U^2-\varepsilon b^2)\big),
\end{eqnarray*}
where
\begin{eqnarray}
\label{R} R= \sqrt{\varepsilon\big((m^2U^2-\varepsilon
b^2)(1-\kappa(m^2U^2-\varepsilon b^2))- m^4 U^2 \dot U^2\big)}.
\end{eqnarray}
So, all surfaces $X_{m,b}$ parametrized by the coordinates
(\ref{coordinaterho})-(\ref{coordinatetheta}) have the same constant
mean curvature $H$ if and only if $U$ satisfies the following
ordinary differential equation
\begin{equation}
\label{H-U} 2H\, R =m^2U \ddot U+\big(m^2 +\frac{\kappa m^4
U^2}{1-\kappa(m^2U^2-\varepsilon b^2)}\big)\, \dot U^2- \varepsilon
\big(1-\kappa(m^2U^2-\varepsilon b^2)\big).
\end{equation}
It is useful now to consider conformal parameters on $\Sigma$ by
changing variables
$$(\tilde s, \tilde \theta)\mapsto (u,v)=: (\int \frac{\dd
\tilde s}{U}, \tilde \theta).$$
Plugging $\partial u/\partial\tilde s=\dd u/\dd \tilde s =1/U,\,
\partial v/\partial\tilde \theta =\dd v/\dd \tilde \theta=1$ into
(\ref{2ndff}) implies that its coefficients are now changed as
\[
\tilde e \mapsto \tilde e U^2, \quad \tilde f\mapsto \tilde f U,
\quad \tilde g \mapsto \tilde g.
\]
The metric induced on $\Sigma$ becomes $U^2(\dd u^2 + \dd v^2)$.
Thus the mean curvature is
\[
2H U^2 = \tilde e U^2 + \tilde g.
\]
So the coefficient $\psi^1$ of the Hopf differential $\Psi^1$ (see
Section 4) in these parameters is written as
\[
\psi^1=\Big(\frac{\tilde e U^2 -\tilde g}{2}\Big)-i\, \tilde f U =
\big(H U^2 -\tilde g\big)-i\, \tilde f U.
\]
Since $\tilde g =-R/m^2$ and $\tilde
f=\frac{b}{m^2U}\sqrt{1-\kappa(m^2U^2-\varepsilon b^2)}=\frac{b}{m^2
U}\cs_\kappa(\rho)$ it follows that
\begin{eqnarray*}
\psi^1=\big(HU^2
+\frac{R}{m^2}\big)-i\,\frac{b}{m^2}\cs_\kappa(\rho).
\end{eqnarray*}
However by the very definition of $R$ the expression
(\ref{integral-final}) reads
\[
\frac{R}{\cs_\kappa (\rho)} = I-2H\int \sn_\kappa (\rho)\, \dd \rho.
\]
So, replacing the identity $\frac{\dd }{\dd
\rho}\,\cs_\kappa(\rho)=-\kappa\, \sn_\kappa(\rho)$ gives
\[
\kappa R = \cs_\kappa(\rho)\,\big(\kappa I +2H\cs_\kappa
(\rho)\big).
\]
We are interested here on $\kappa\neq 0$ (it is a well-known fact
that $\psi^1$ is holomorphic for $\kappa=0$). In this case it holds
that
\begin{eqnarray*}
& & \kappa \psi^1 = \big( \kappa H U^2
+\frac{1}{m^2}\cs_\kappa(\rho)\,\big(\kappa I +2H\cs_\kappa
(\rho)\big)\big)-i\,\frac{\kappa b}{m^2}\, \cs_\kappa (\rho)\\
& & =\frac{1}{m^2}\,\big(\kappa H \sn_\kappa^2(\rho)+\kappa
H\varepsilon b^2+\kappa I \cs_\kappa(\rho) + 2H \cs_\kappa^2(\rho)
\big)-i\,\frac{\kappa b}{m^2}\, \cs_\kappa (\rho).
\end{eqnarray*}
Now we want to compute  the coefficient $\psi^2$  of the
differential $\Psi^2$ on the conformal coordinates $u,v$ defined
just above. We have
\[
\partial_u = X_{\tilde s}\frac{\partial \tilde s}{\partial u}=X_{\tilde
s} U.
\]
Using equations (\ref{X}) and (\ref{tildes}) one proves that
\[
\langle\partial_u,\partial_t\rangle =\varepsilon\frac{\dd t}{\dd
\tilde s} U\Big(\frac{m^2 U^2 -\varepsilon b^2}{m^2 U^2}\Big).
\]
However
\begin{eqnarray*}
& & \big(m^2 U^2 -\varepsilon b^2\big)\frac{\dd t}{\dd \tilde
s}=\Big(\frac{\sn_\kappa^2(\rho)}{mU}\frac{\dd t}{\dd \tilde
s}\Big)\, mU=\Big(\frac{\sn_\kappa^2(\rho)}{W}\frac{\dd t}{\dd
\rho}\Big)\, mU\\
& & =mU\big(I-2H\int\sn_\kappa(\rho)\, \dd \rho\big).
\end{eqnarray*}
Therefore
\begin{eqnarray*}
& & \kappa \big(m^2 U^2 -\varepsilon b^2\big)\frac{\dd t}{\dd \tilde
s}=mU\big(\kappa I+ 2H\cs_\kappa(\rho)\big).
\end{eqnarray*}
We conclude that
\begin{eqnarray*}
& & \kappa \langle\partial_u,\partial_t\rangle
=\frac{\varepsilon}{m}\,\big(\kappa I+ 2H\cs_\kappa(\rho)\big).
\end{eqnarray*}
We also compute
\begin{eqnarray*}
& & \partial_v =X_{\tilde\theta}=X_\theta \frac{\partial
\theta}{\partial\tilde\theta}= \frac{1}{m}\Big(\partial_\theta +
b\partial_t\Big)
\end{eqnarray*}
and
\[
\langle \partial_v,\partial_t\rangle = \frac{\varepsilon b}{m}.
\]
Thus it results that
\[
\kappa^2 \langle
\partial_u,\partial_t\rangle^2-\kappa^2\langle\partial_v,\partial_t\rangle^2
=\frac{1}{m^2}\,\big(\kappa I+
2H\cs_\kappa(\rho)\big)^2-\frac{\kappa^2 b^2}{m^2}
\]
and
\[
\kappa \langle\partial_u,\partial_t\rangle\langle
\partial_v,\partial_v\rangle=\frac{b}{m^2}\,\big(\kappa I+
2H\cs_\kappa(\rho)\big).
\]
Now since that $\frac{\epsilon}{r^2}=\kappa$ we write
\begin{eqnarray*}
& & \frac{\epsilon}{r}\,\psi^2 =\frac{1}{2}\Big(\kappa\varepsilon
\langle
\partial_u,\partial_t\rangle^2 -\kappa\varepsilon
\langle \partial_v,\partial_t\rangle^2\Big)-i\, \kappa\varepsilon
\langle \partial_u,\partial t\rangle\langle
\partial_v,\partial_t\rangle\\
& & =\frac{\varepsilon}{2\kappa m^2}\Big(\kappa^2 I^2+
4H^2\cs^2_\kappa(\rho)+4H\kappa I \cs_\kappa(\rho)-\kappa^2
b^2\Big)-i\,\frac{\varepsilon b}{m^2}\big(\kappa I+
2H\cs_\kappa(\rho)\big).
\end{eqnarray*}
Therefore
\begin{eqnarray*}
& & 2H\psi^1 -\varepsilon\frac{\epsilon}{r}\, \psi^2
=\frac{1}{m^2}\big(2H^2(\frac{1}{\kappa}+ \varepsilon b^2
)+\frac{1}{2} \kappa (b^2 -I^2)\big) +i\,\frac{b\kappa I}{m^2}.
\end{eqnarray*}
For $\kappa=0$ we have $\cs_\kappa(\rho)=1$ and
(\ref{integral-final}) becomes $R = I-H^2 \rho^2.$ So
\begin{eqnarray*}
& & \psi^1 = \frac{1}{m^2}\big(Hm^2U^2 +R\big)-i\,\frac{b}{m^2}
=\frac{1}{m^2}\big(H\varepsilon b^2 +I \big)-i\,\frac{b}{m^2}.
\end{eqnarray*}
Thus, the differential $Q$ has constant coefficient for any surface
on the family $X_{m,b}:\Sigma\to \Mm\times\Rr$  of screw-motion
invariant CMC surfaces on $\Mm\times \Rr$  starting  (for $m=1$)
from some given CMC surface. Its final expression is:
\[
\psi = -\frac{1}{2m^2\kappa}\big(\kappa^2 I^2-4H^2-\kappa
b^2(4H^2\varepsilon+\kappa) \big)+\,i\frac{b\kappa I}{m^2}
\]
for $\kappa\neq 0$. From the same calculations, we assure that the
Hopf differential has constant coefficient for $\kappa=0$:
\[
\psi^1=\frac{1}{m^2}\big(H\varepsilon b^2 +I \big)-i\,\frac{b}{m^2}
\]
In the case $\kappa\neq 0$, we have for rotational examples ($b=0$,
$m=1$) that
\[
\psi =-\frac{1}{2\kappa}\big(\kappa^2 I^2-4H^2 \big)
\]
Thus $Q=0$ for CMC rotational examples if and only if $4H^2-\kappa^2
I^2=0$ or
\[
I = \pm\frac{2H}{\kappa}.
\]
We now determine explicitly the CMC rotational examples with $Q=0$.
In order to do this, we replace  $I=\pm 2H/\kappa$  in
(\ref{integral}). Since $W^2 =\sn_\kappa^2(\rho)(1+\varepsilon
\frac{\dd t}{\dd \rho}^2)$ it follows that
\[ \frac{\frac{\dd t}{\dd \rho}}{\sqrt{1+\varepsilon \frac{\dd t}{\dd\rho}^2}}\,\sn_\kappa(\rho)=
-\frac{2H}{\kappa}\,\big(\pm 1 +\kappa\int \, \sn_\kappa (\rho)\,\dd
\rho\big)=-\frac{2H}{\kappa}\big(\pm 1-\cs_\kappa(\rho)\big).
\]
So squaring both sides and taking inverses
\[
\frac{1+\varepsilon \frac{\dd t}{\dd\rho}^2}{\big(\frac{\dd t}{\dd
\rho}\big)^2}=\frac{\kappa}{4H^2}\,
\frac{\kappa\sn^2_\kappa(\rho)}{\big(\pm
1-\cs_\kappa(\rho)\big)^2}=\frac{\kappa}{4H^2}\,
\frac{1-\cs^2_\kappa(\rho)}{\big(\pm 1-\cs_\kappa(\rho)\big)^2}.
\]
Thus for  $I=-2H/\kappa$ one has
\[
\big(\frac{\dd \rho}{\dd t}\big)^2+\varepsilon
=\frac{\kappa}{4H^2}\,
\frac{1+\cs_\kappa(\rho)}{1-\cs_\kappa(\rho)}.
\]
However
\begin{eqnarray*}
& &
\frac{1+\cs_\kappa(\rho)}{1-\cs_\kappa(\rho)}=\frac{1}{\kappa}\,\ct^2_\kappa
(\rho/2)=\frac{1}{\kappa}\, \frac{1}{r^2}.
\end{eqnarray*}
Here $\textrm{ct}_\kappa(\rho) = \dot{ \textrm{sn}}_\kappa(\rho)/
\textrm{sn}_\kappa(\rho)$ is the geodesic curvature of the geodesic
circle centered at $p_0$ with radius $\rho$ in
$\mathbb{M}^2(\kappa)$ and $r$ is the Euclidean radial distance $r$
measured from $p_0$ on the Euclidean model for $\Mm$. Thus for
$I=-2H/\kappa$ we have
\begin{eqnarray*}
& & \big(\frac{\dd \rho}{\dd t}\big)^2+\varepsilon =
\frac{1}{4H^2}\,\frac{1}{r^2}.
\end{eqnarray*}
Now
\[
\big(\frac{\dd \rho}{\dd t}\big)^2=\big(\frac{\dd \rho}{\dd
r}\big)^2\,\big(\frac{\dd r}{\dd t}\big)^2=\frac{4}{(1+\kappa
r^2)^2}\,\big(\frac{\dd r}{\dd t}\big)^2.
\]
So the resulting equation is
\[
\frac{2}{1+\kappa r^2}\,\frac{2Hr\,\dd r}{\sqrt{1- 4H^2 \varepsilon
r^2}}=\dd t.
\]
We change variables defining $1+\kappa r^2 =u$. We then change
variables  again by defining ($\kappa<0$) $v=\varepsilon u
-\big(\varepsilon+\kappa/4H^2\big)$ and
$v=\big(\varepsilon+\kappa/4H^2\big)-\varepsilon u$ (for
$\kappa>0$). Next, we put $w=\sqrt v$. So $\dd v =2w\,\dd w$ and the
final form of the equation is
\[
\frac{2\dd
w}{w^2+\big(\varepsilon+\kappa/4H^2\big)}=-\sqrt{-\kappa}\,\dd t,\,
(\kappa<0),\quad \frac{2\dd
w}{w^2-\big(\varepsilon+\kappa/4H^2\big)}=\sqrt{\kappa}\,\dd t, \,
(\kappa>0).
\]
We suppose that $4H^2\varepsilon +\kappa>0$. Then writing $c^2 =
\varepsilon+\kappa/4H^2$ one has
\[
\frac{2}{c}\, \arctan (w/c)= -\sqrt{-\kappa}\, t,\quad \kappa<0
\]
and
\[
\frac{1}{c}\log\Big|\frac{w+c}{w-c}\Big|=-\sqrt{\kappa}\, t,\quad
\kappa>0
\]
In this last case, notice that $|w|<c$ (respectively, $|w|>c$) if
and only if $\varepsilon=1$ (resp., $\varepsilon=-1$). We fix
initially $\kappa<0$. Then necessarily $\varepsilon=1$ and
\[
v=w^2=-c^2\kappa\, \ct_{-\kappa}^{-2}(-ct/2)
\]
so that
\[
\varepsilon u = v+c^2 = c^2\,\big(1-\kappa\,
\ct_{-\kappa}^{-2}(-ct/2)\big)=c^2\cs_{-\kappa}^{-2}(ct/2).
\]
Since $u=1+\kappa r^2 =1/\cs_\kappa^2(\rho/2)$ and   $\varepsilon
c^2 = 1+\kappa \varepsilon/4H^2=1+\kappa/4H^2\varepsilon$ then
\begin{equation}
\label{1+-} \big(4H^2\varepsilon+\kappa\big)\,
\sn_\kappa^2(\rho/2)+4H^2\varepsilon\, \sn_{-\kappa}^2
(ct/2)=1,\quad (\kappa<0),
\end{equation}
where $c=\sqrt{\varepsilon+\kappa/4H^2}$ and $\varepsilon=1$. The
same formula holds for $\kappa>0$, $\varepsilon=1$. We have for
$\kappa>0,\, \varepsilon=-1$ that
\begin{equation}
\label{-1+-} 4H^2\varepsilon \kappa
\,\sn_{-\kappa}^2(ct/2)=-\big(4H^2\varepsilon+\kappa\big)\,
\cs_\kappa^2(\rho/2).
\end{equation}
We now treat the case $\varepsilon +\kappa/4H^2<0$. We denote
$c^2=-(\varepsilon +\kappa/4H^2)$. Thus for $\kappa>0$ and
$\varepsilon=-1$ the solution is
\begin{equation}
\label{-1--} \big(4H^2\varepsilon+\kappa\big)\,
\sn_\kappa^2(\rho/2)-4H^2\varepsilon \sn_\kappa^2(ct/2)=1
\end{equation}
The same formula holds for $\kappa<0$, $\varepsilon=-1$ when we have
$|w|<c$. For $\varepsilon=1$, we necessarily have $\kappa<0$ and
$|w|>c$. Thus
\begin{equation}
\label{1--} 4H^2\varepsilon \kappa
\sn_\kappa^2(ct/2)=\big(4H^2\varepsilon+\kappa\big)\, \cs_\kappa^2
(\rho/2).
\end{equation}
Finally for $\varepsilon +\kappa/4H^2=0$ one obtains
\begin{equation}
\label{0} t^2 = \epsilon\frac{4}{\kappa}\, \cs_\kappa^2(\rho/2).
\end{equation}
Next, we consider $I=2H/\kappa$. For this choice we have
\[
\big(\frac{\dd \rho}{\dd t}\big)^2+\varepsilon
=\frac{\kappa}{4H^2}\,
\frac{1-\cs_\kappa(\rho)}{1+\cs_\kappa(\rho)}.
\]
So the resulting equation  is
\[
\frac{2}{1+\kappa r^2}\,\frac{2H\,\dd r}{\sqrt{\kappa^2 r^2- 4H^2
\varepsilon }}=\dd t
\]
We change variables considering $u=1/r^2+\kappa =
1/\sn_\kappa^2(\rho/2)$. We then change variables again by defining
$v=\kappa\big(\varepsilon +\kappa/4H^2\big)-\varepsilon u$. Finally
we put $w^2=v$. So
\[
\frac{2\dd w}{\kappa\big(\varepsilon +\kappa/4H^2\big)-w^2}=\dd t.
\]
First, we consider $c^2=\varepsilon +\kappa/4H^2>0$. In this case
there are no examples with $\kappa<0$. For $\kappa>0$ and
$\varepsilon=1$
\begin{equation}
\label{1++} \big(4H^2\varepsilon+\kappa\big)\,\kappa
\sn_\kappa^2(\rho/2)=4H^2\varepsilon \,
\cs_{-\kappa}^2(\sqrt{\varepsilon+\kappa/4H^2}\, t/2).
\end{equation}
For $\kappa>0$ and $\varepsilon=-1$
\begin{equation}
\label{-1++} \big(4H^2\varepsilon+\kappa\big)\,
\sn_\kappa^2(\rho/2)=-4H^2\varepsilon \,
\sn_{-\kappa}^2(\sqrt{\varepsilon+\kappa/4H^2}\, t/2).
\end{equation}
Now, we consider the case  $-c^2=\varepsilon +\kappa/4H^2<0$. For
$\kappa<0$ and $\varepsilon=1$ we have
\begin{equation}
\label{1-+}
\big(4H^2\varepsilon+\kappa\big)\kappa\,\sn_\kappa^2(\rho/2)=4H^2\varepsilon
\cs_{\kappa}^2\big(\sqrt{-\big(\varepsilon +\kappa/4H^2\big)}
t/2\big)
\end{equation}
The same expression holds for $\kappa>0,\, \varepsilon=-1$. For
$\kappa<0,\,\varepsilon=-1$ we have
\begin{equation}
\label{-1-+}
\big(4H^2\varepsilon+\kappa\big)\,\sn_\kappa^2(\rho/2)=4H^2\varepsilon
\sn_{\kappa}^2\big(\sqrt{-\big(\varepsilon +\kappa/4H^2\big)}
t/2\big)
\end{equation}

\vspace{0.3cm}

\noindent {\bf Theorem 2.} {\it The revolution  surfaces with
constant mean curvature $H$ and $Q=0$ on $\Mm\times \Rr$ correspond
to the values $I=-2H/\kappa,\,2H/\kappa$. These surfaces are
described by the formulae (\ref{1+-})-(\ref{-1-+}) just above.}

\vspace{0.3cm}

For $\varepsilon=1$, the formulae above were already obtained in
\cite{AR}  by other integration methods.

\subsection{Solving the mean curvature equation}

We proved on Section 2.2 that a given  CMC helicoidal surface could
be deformed on isometric   helicoidal surfaces with the same mean
curvature. In this section, we give explicit parameterizations to
these families.

We denote in what follows the variable $\tilde s$ simply as $s$.
Squaring both sides of (\ref{integral-final}) one finds
\begin{eqnarray}
\label{integral3} \varepsilon\,\frac{\big(m^2U^2-\varepsilon
b^2\big)\big(1-\kappa(m^2U^2-\varepsilon b^2)\big)-m^4U^2\dot U^2}
{1-\kappa(m^2U^2-\varepsilon b^2)}= (2H\int \sn_\kappa(\rho)\, \dd
\rho -I)^2.
\end{eqnarray}
In particular, for $\kappa=0$ since $\sn_\kappa(\rho)=\rho$ and
$\rho^2 = m^2 U^2-\varepsilon b^2$ then (\ref{integral3}) becomes
\begin{equation}
\label{kappa0} \varepsilon\,\big(m^2U^2-\varepsilon b^2-m^4U^2\dot
U^2\big) = (Hm^2 U^2-H\varepsilon b^2 -I)^2.
\end{equation}
For $\varepsilon=1$ after the substitutions $x=: mU$ and $z=:x^2$
this equation reads
\begin{equation}
\label{eq-euclidean} \frac{\dot z^2}{4}= -H^2z^2+(1+2Ha)z-(a^2+b^2),
\end{equation}
where $a= H\varepsilon b^2 +I$. This equation was solved in
\cite{dCD} and its solutions completely integrated. For
$\varepsilon=-1$ the same substitutions show that (\ref{kappa0})
becomes
\begin{equation} \label{eq-lorentzian} \frac{\dot z^2}{4}=
H^2z^2+(1-2Ha)z +(a^2+b^2).
\end{equation}
Completing squares this equation reads
\begin{eqnarray}
& & \frac{\dot
z^2}{4H^2}=\Big(z+\frac{1-2Ha}{2H^2}\Big)^2+\frac{4H^2b^2+4Ha-1}{4H^4},
\end{eqnarray}
for $H\neq 0$ and $\dot z^ 2/4 = z +(a^2+b^2)$ for $H=0$.  This last
equation may be rewritten as
\[
\frac{\dd z}{\sqrt{z+(a^2+b^2)}} = 2\dd s
\]
whose solution is of the form $m^2 U^2 = z = s^2-(a^2+b^2)$, where
$a=I$ since $H=0$. This family contains a Lorentzian catenoid as
initial surface. In fact, considering the values $m=1$ and $b=0$, we
have $U^2= s^2-I^2$ and $\rho^2=U^2$. So $s=\sqrt{\rho^2+I^2}$ and
$\dd s=(\rho/\sqrt{\rho^2+I^2})\,\dd\rho$. The expression
(\ref{coordinatet}) reads
\begin{equation*}
t  =\int\frac{I}{\rho}\, \dd s=\int\frac{I}{\sqrt{\rho^2+I^2}}\, \dd
\rho =I \, \textrm{arcsinh}\big(\rho/I\big)
\end{equation*}
Thus the (half of the) catenoid is described as the graph of
\begin{equation}
\label{catenoidt} \rho=\rho(t)=I\, \sinh \big(t/I\big)
\end{equation}
We remark that this curve is singular at $t=0$ and asymptotes a
light cone there. For the catenoid we have $\tilde \theta =\theta$.
We now describe the family associated to such a catenoid  by the
integrals (\ref{coordinaterho})-(\ref{coordinatetheta}). For the
other members of the family that evolves from the Lorentzian
catenoid we have $\rho^2=m^2U^2+b^2=s^2-(I^2+b^2)+b^2 = s^2-I^2$ and
$s=\sqrt{\rho^2+I^2}$. So
\begin{equation}
\label{catenoidmt} t=\int \sqrt{\frac{\rho^2+b^2}{\rho^2+I^2}}\,
\frac{I}{\rho}\, \dd \rho
\end{equation}
and the coordinate $\tilde \theta (\rho, \theta)$ is given by
\begin{equation}
\label{catenoidmtheta} \tilde\theta(\rho, \theta) = m \,\theta
-mbI\int \frac{1}{\rho^2\,\sqrt{\rho^2 -I^2}\sqrt{\rho^2+b^2}}\, \dd
\rho.
\end{equation}

Turning back to the Lorentzian equation (\ref{eq-lorentzian}) for
$H\neq 0$, if we
consider $w=z+\frac{1-2Ha}{2H^2}$ and
$c^2=|\frac{4H^2b^2+4Ha-1}{4H^4}|$ we have
\[
\int \frac{\dd w}{\sqrt{w^2\pm c^2}}= 2Hs
\]
whose general solutions are, for sign $+$
\begin{equation}
\label{HL3+} m^2U^2 =c\,\sinh\big(2H(s-s_0)\big)+\frac{1-2Ha}{2H^2}
\end{equation}
and for sign $-$
\begin{equation}
\label{HL3-} m^2 U^2
=c\,\cosh\big(2H(s-s_0)\big)+\frac{1-2Ha}{2H^2},
\end{equation}
where
\[
c=\big|\frac{4H^2b^2+4Ha-1}{4H^4}\big|^{1/2}.
\]
We may make explicit the parametrization describing both $U$ and
$U\dot U$ in terms of these solutions.

\vspace{0.3cm} \noindent {\bf Theorem 3. }{\it A family of maximal
space-like helicoidal surfaces in $\mathbb{L}^3$ containing a
Lorentzian catenoid is described by the formulae
(\ref{catenoidt})-(\ref{catenoidmtheta}). The formulae (\ref{HL3+})
and (\ref{HL3-}) describe families of helicoidal CMC surfaces on
$\mathbb{L}^3$.} \vspace{0.3cm}

Now, we consider the case $\kappa\neq 0$. Since  $\frac{\dd }{\dd
\rho}\,\cs_\kappa(\rho)=-\kappa\, \sn_\kappa(\rho)$ then
\begin{equation}
\label{rhs} (-2H\kappa\int \sn_\kappa(\rho)\, \dd \rho+\kappa I)^2=
4H^2 \cs_\kappa^2(\rho)+4H \kappa I\cs_\kappa(\rho)+\kappa^2 I^2.
\end{equation}
Since $\cs_\kappa(\rho) = (1-\kappa(m^2U^2-\varepsilon b^2))^{1/2}$,
defining $z=:(1-\kappa(m^2U^2-\varepsilon b^2))^{1/2}$ for
$\kappa\neq 0$ one finds $z^2-1 -\varepsilon\kappa b^2 = -\kappa m^2
U^2$. Therefore $z \dot z= -\kappa m^2 U \dot U$ which implies that
$\kappa^2 m^4 U^2 \dot U^2 = z^2 \dot z^2.$  Multiplying both sides
of the expression (\ref{integral3}) by $\kappa^2$ and replacing the
expression (\ref{rhs}) on the right hand side of the resulting
equation we obtain a first integral to the equation (\ref{H-U})
\begin{eqnarray}
\label{integral4} \nonumber & &
\kappa^2\varepsilon\,\frac{\big(m^2U^2-\varepsilon
b^2\big)\big(1-\kappa(m^2U^2-\varepsilon b^2)\big)-m^4U^2\dot U^2}
{1-\kappa(m^2U^2-\varepsilon b^2)}= \big(2H
\big(1-\kappa(m^2U^2\\
& &-\varepsilon b^2)\big)^{1/2} +\kappa I\big)^2.
\end{eqnarray}
In terms of $z$ this equation reads
\begin{eqnarray}
\label{kappa} \dot z^2 = -(4H^2\varepsilon+\kappa)\, z^2 - 4H \kappa
I \varepsilon\,z +\kappa(1-\kappa I^2\varepsilon).
\end{eqnarray}
If we assume that $4H^2\varepsilon+\kappa\neq 0$ then we obtain
after completing squares that
\begin{eqnarray}
\frac{\dot z^2}{4H^2\varepsilon +\kappa} =-\Big(z+\frac{2H\kappa I
\varepsilon}{4H^2\varepsilon+\kappa}\Big)^2
+\frac{\kappa}{(4H^2\varepsilon+\kappa)^2}\, \big(4H^2\varepsilon
+\kappa -\kappa^2 I^2\varepsilon\big).
\end{eqnarray}
We first consider the case $4H^2\varepsilon+\kappa<0$. If
$\kappa\big(4H^2\varepsilon +\kappa -\kappa^2 I^2\varepsilon\big)<0$
then putting $w=z+\frac{2H\kappa I
\varepsilon}{4H^2\varepsilon+\kappa}$ we get
\begin{eqnarray}
-\frac{\dot w^2}{4H^2\varepsilon+\kappa}=w^2+c^2,
\end{eqnarray}
where
\[
c^2=\big|\frac{\kappa}{(4H^2\varepsilon+\kappa)^2}\,
\big(4H^2\varepsilon +\kappa -\kappa^2 I^2\varepsilon\big)\big|.
\]
The general solution is in this case
\begin{eqnarray}
\label{--} z = c \sinh\big((-4H^2\varepsilon-\kappa)^{1/2}
(s-s_0)\big)-\frac{2H\kappa I \varepsilon}{4H^2\varepsilon+\kappa}.
\end{eqnarray}
If $\kappa\big(4H^2\varepsilon +\kappa -\kappa^2
I^2\varepsilon\big)>0$  then
\begin{equation*}
 -\frac{\dot w^2}{4H^2\varepsilon+\kappa}=w^2 -c^2
\end{equation*}
with solution given by
\begin{equation}
\label{-+} z = c \cosh\big((-4H^2\varepsilon-\kappa)^{1/2}
(s-s_0)\big)-\frac{2H\kappa I \varepsilon}{4H^2\varepsilon+\kappa}.
\end{equation}
Now we consider the case $4H^2\varepsilon+\kappa>0$. Here  we
necessarily have $\kappa\big(4H^2\varepsilon +\kappa -\kappa^2
I^2\varepsilon\big)>0$. The equation becomes
\[
\frac{\dot w^2}{4H^2\varepsilon+\kappa}=c^2-w^2,
\]
whose solution is
\begin{equation}
\label{++} z = c
\sin\big((4H^2\varepsilon+\kappa)^{1/2}(s-s_0)\big)-\frac{2H\kappa I
\varepsilon}{4H^2\varepsilon+\kappa}.
\end{equation}
It remains to see what happens for $4H^2\varepsilon+\kappa=0$. In
this case the equation becomes
\[
\dot z^2 =- 4H \kappa I \varepsilon\,z +\kappa(1-\kappa
I^2\varepsilon).
\]
If $H\kappa I=0$ then we have necessarily $\kappa(1-\kappa
I^2\varepsilon)>0$ and
\begin{equation}
\label{00} z =\big(\kappa(1-\kappa
I^2\varepsilon)\big)^{1/2}(s-s_0).
\end{equation}
When $H\kappa I\neq 0$ then the equation is
\begin{equation*}
\frac{\dd z}{\sqrt{- 4H \kappa I \varepsilon\,z +\kappa(1-\kappa
I^2\varepsilon)}}=\dd s
\end{equation*}
with solution
\begin{equation}
\label{0n0} z=-\frac{1}{4H\kappa I \varepsilon}\, \big(\frac{1}{4}(s
-s_0)^2 -\kappa (1-\kappa I^2\varepsilon)\big).
\end{equation}

\vspace{0.3cm} \noindent{\bf Theorem 4.} {\it The formulae
(\ref{--})-(\ref{0n0}) describe two-parameter families of helicoidal
CMC examples on $\Mm\times\Rr$.}

\vspace{0.3cm}

For $\varepsilon=1$, the formulae above were previously obtained in
\cite{ET}.

\section{Rotationally invariant CMC discs on Lorentzian products}

\subsection{Qualitative description}

In this section we  consider only space-like revolution surfaces in
Lorentzian products $\Mm\times \Rr$. We assume that the parameter
$s$ on (\ref{X}) is the arc length of the profile curve. So,
$\dot\rho^2 -\dot t^2 =1$. We denote by $\varphi$ the hyperbolic
angle with the horizontal axis $\partial_\rho$. So, $\Sigma$ has
constant mean curvature $H$ if and only if
$(\rho(s),t(s),\varphi(s))$ is solution to  the following ordinary
differential equations system
\begin{eqnarray}
\label{system}
& &\dot\rho = \cosh \varphi,\nonumber\\
& & \dot t = \sinh\varphi,\nonumber\\
& & \dot\varphi = - 2H -\sinh\varphi\, \textrm{ct}_\kappa(\rho).
\end{eqnarray}
The {\it flux} $I'$ through an horizontal plane $\Pp
=\mathbb{M}^2(\kappa)\times \{t\}$ is, up to a constant, given by
the expression for $I$ in terms of $s$:
\begin{equation}
\label{flux}I' = I+\frac{2H}{\kappa} = \sinh
\varphi\,\sn_\kappa(\rho) + 2H\, \int^\rho_0 \sn_\kappa(\tau)\,
\dd\tau.
\end{equation}
Integrating the last term on (\ref{flux}) one obtains
\begin{eqnarray}
\label{flux2} & & I'= \sinh\varphi \,\sn_\kappa (\rho) + 4H
\sn_\kappa^2 (\frac{\rho}{2}).
\end{eqnarray}
The solutions  for (\ref{system}) for which $Q=0$ vanishes are those
with $I=\pm \frac{2H}{\kappa}$ or  $I'=0, \frac{4H}{\kappa}$. We
give later a qualitative description  of these solutions.

Since that $\cosh \varphi$ never vanishes on the maximal interval
for a solution to (\ref{system}) it follows that
\[
\frac{\dd t}{\dd \rho} = \frac{\dd t}{\dd s}\,\frac{\dd s}{\dd \rho}
= \tanh \varphi.
\]
Denoting $u = \sinh \varphi$ we also obtain
\[
\frac{\dd u}{\dd \rho} = \frac{\dd u}{\dd \varphi}\,\frac{\dd
\varphi}{\dd s}\,\frac{\dd s}{\dd \rho} =  -2H -\sinh\varphi\,
\ct_\kappa (\rho).
\]
Thus the system (\ref{system}) above is equivalent to
\begin{eqnarray}
\label{system2} & & \frac{\dd t }{\dd \rho} = \frac{u}{\sqrt{ 1 +
u^2}},\nonumber\\
& & \frac{\dd u}{\dd \rho} = -2H - u\, \ct_\kappa(\rho).
\end{eqnarray}
It is clear that solutions to the system (\ref{system2}) are defined
on the whole real line and the profile curve may be written as a
graph over the $\rho$-axis. Now, we  begin describing  the maximal
solutions, i.e., solutions for $H=0$. If we consider a fixed value
for $I'$ then the condition $H=0$ implies that
\begin{equation}
\label{maximal} I' = \sinh \varphi \, \sn_\kappa (\rho)
\end{equation}
So, the horizontal planes are the unique maximal revolution surfaces
with $I'=0$. In fact if we put $I'=0$ at (\ref{maximal}) we have
$\sinh\varphi=0$ for $\rho>0$. Thus, $\dot t =0$ and we conclude
that the solution is an horizontal plane. Hence, we may assume
$I\neq 0$. In this case, since that $\sn_\kappa (\rho)\to 0$ if
$\rho\to 0$ is follows that $\sinh \varphi \to \infty$ if $\rho\to
0$. So, $\Sigma$ has a singularity and asymptotes the light cone at
$p_0$ (the light cone corresponds to $\varphi =\infty$). Moreover
$\sinh\varphi \to 0$ if $\rho\to \infty$ in the case $\kappa\le 0$.
This means that these maximal surfaces asymptotes an horizontal
plane for $\rho\to \infty$, i.e., these surfaces have planar ends.
These examples are not complete in the spherical case $\kappa>0$,
since we have $\sinh\varphi\to \infty$ if $\rho\to \frac{\pi}{\sqrt
\kappa}$.

Consider now the case $H\neq 0$. We observe that the solutions for
(\ref{system2}) have no positive minimum  for $\rho$. Otherwise, the
solutions must have vertical tangent plane at the minimum points
(this is impossible since the solutions are space-like and, in fact,
are graphs over the horizontal axis). Hence, the unique possibility
for the existence of a isolated singularity  is that $\rho \to 0$.
In this case the solutions are regular if and only if the
$\varphi\to 0$ as $\rho \to 0$ what implies that  $\sinh \varphi \to
0$ as $\rho \to 0$. So, necessarily $I'=0$ as we could see taking
the limit $\rho \to 0$ in (\ref{flux2}) above. So, examples of
solutions for the systems above which touch orthogonally the
revolution axis have $I'=0$. Reciprocally, if we put $I'=0$ in
(\ref{flux2}) we get
\begin{eqnarray*}
& & 0= \sinh\varphi \,\sn_\kappa (\rho) + 4H \sn_\kappa^2
(\frac{\rho}{2}).
\end{eqnarray*}
So, dividing the expression above by $2\sn^2_\kappa
(\frac{\rho}{2})$ we have
\begin{equation}
\label{fluxzero} \sinh\varphi\,\ct_\kappa (\frac{\rho}{2}) = -2H.
\end{equation}
One easily verifies that $\sinh \varphi \to 0$ if $\rho\to 0$. So
all solutions for (\ref{system2}) with $I'=0$ reach the revolution
axis orthogonally as we noticed earlier. Thus these solutions
correspond to initial conditions $t(0) = t_0$, $\rho(0) =0$ and
$\varphi(0) = 0$ for the system (\ref{system}). Now we have
\begin{eqnarray*}
\ct_\kappa (\rho)  = \frac{1}{2}\big(-\frac{2H}{\sinh\varphi}+
\kappa\, \frac{\sinh\varphi}{2H}\big) = \frac{-4H^2 +
\kappa\,\sinh^2\varphi}{4H \sinh\varphi}.
\end{eqnarray*}
Replacing this on the third equation on (\ref{system}) we obtain
\begin{equation}
\label{asymptoticangle} \frac{\dd \varphi}{\dd s} = \frac{1}{4H}
\,(-4H^2 -\kappa\, \sinh^2\varphi).
\end{equation}
We observe that $\dot\varphi = -H$ is the corresponding equation for
the case $\kappa =0$, i.e., for hyperbolic spaces in $\mathbb{L}^3$.
This could be obtained as a limiting case if we take $\kappa \to 0$.
For $\kappa<0$, the range for the angle $\varphi$ is $0\le \varphi <
\varphi_\infty = \textrm{arcsinh}(2|H|/\sqrt{-\kappa})$. The surface
necessarily asymptotes a spacelike cone with angle $\varphi_\infty$.
Indeed the equation (\ref{asymptoticangle}) is equivalent to
\[
\frac{1}{4H}\int_0^{\varphi_\infty} \frac{\dd \varphi}{-4H^2
-\kappa\,\sinh^2 \varphi} = \int_0^{\infty} \dd s = \infty.
\]
There are no complete solutions for $\kappa>0$ and $H\neq 0$, since
that the angle at $\rho=0$ and at $\rho= \frac{\pi}{\sqrt \kappa}$
are not the same unless we have $H=0$.

Finally, we study the case when $\varphi \to \varphi_0$ as $\rho \to
0$ for some positive value of $\varphi_0$. This means that the
solution asymptotes a space-like cone at $p_0$.  In this case $\sinh
\varphi \to \sinh \varphi_0 < \infty$ as $\rho\to 0$. Thus taking
the limit $\rho \to 0$ in (\ref{flux2}) we obtain $I'=0$. So, as we
seen above, necessarily $\varphi_0 =0$. This contradiction implies
that there are no examples with $\varphi_0>0$.

It remains to give a look at the case $\varphi \to \infty$ as $\rho
\to 0$. In this case, the solution asymptotes the light cone at
$p_0$. For any non zero value of $I'$, we obtain after dividing
(\ref{flux2}) by $\sn_\kappa^2 (\rho/2)$ and taking limit for $\rho
\to \infty$ that $\sinh \varphi \to 2|H|/\sqrt{-\kappa}$. Moreover,
the angle $\varphi$ is always decreasing in the range
$(2|H|/\sqrt{-\kappa}, \infty)$ as $\rho$ increases in $(0,
\infty)$. For example, consider the values $\kappa< 0$ and $I'
=\frac{4H}{\kappa}$. Replacing this value for $I'$ in (\ref{flux})
we get
\[
0=\sinh\varphi \,\sn_\kappa(\rho) + 4H \big(\sn_\kappa^2
(\frac{\rho}{2})-\frac{ 1}{\kappa}\big).
\]
So we conclude that
\begin{equation}
\label{fluxH} \kappa\,\sinh\varphi = 2H \,\ct_\kappa
\big(\frac{\rho}{2}\big).
\end{equation}
Thus the solution satisfies $\sinh\varphi\to \infty$ if $\rho\to 0$.
This means that $\Sigma$ asymptotes the light cone at the point
$p_0$. Moreover, we have that $\sinh\varphi\to 2|H|/\sqrt{-\kappa}$
if $\rho \to \infty$. Replacing (\ref{fluxH}) at the third equation
in (\ref{system}) we obtain
\begin{eqnarray*}
\ct_\kappa (\rho) = \frac{1}{2} \big(\kappa\,\frac{\sinh\varphi}{2H}
- \frac{2H}{\sinh\varphi}\big) = \frac{-4H^2 +
\kappa\,\sinh^2\varphi}{4H \sinh\varphi}
\end{eqnarray*}
and
\[
\frac{\dd \varphi}{\dd s} = \frac{1}{4H}\,\big(-4H^2 - \kappa
\sinh^2 \varphi\big).
\]
Since $\varphi$ satisfies $\sinh \varphi >\sinh \varphi_\infty =
\frac{2|H|}{\sqrt{-\kappa}}$ then we conclude that $\dot \varphi <0$
for all $s$. So, the angle decreases from $\infty$ at $\rho \to 0$
to its infimum value $\varphi_\infty$ as $\rho \to \infty$.

We summarize the facts above in the following theorem.

\vspace{0.3cm}

\noindent {\bf Theorem 5.} {\it Let $\Sigma$ be a rotationally
invariant surface with constant mean curvature $H$ in the Lorentzian
product $\mathbb{M}^2(\kappa) \times \Rr$ with $\kappa\le 0$. If
$H=0$ either $\Sigma$ is a horizontal plane $\Pp=\mathbb{M}^2
(\kappa) \times \{t\}$ or $\Sigma$ asymptotes a light cone with
vertex at some point $p_0$ of the rotation axis.  In this  case,
$\Sigma$ has a singularity at $p_0$ and has horizontal planar ends.
We refer to these singular surfaces as {\rm Lorentzian catenoids}.

If $H\neq 0$ either $\Sigma$ is a complete disc-type surface meeting
orthogonally the rotation axis or $\Sigma$ asymptotes a light cone
with vertex $p_0$ at the rotation axis. In the first case, the angle
between the surface and the horizontal planes asymptotes
$2|H|/\sqrt{-\kappa}$ as the surface goes to the asymptotic boundary
$\partial_\infty \Mm\times \Rr$. In the last case, the surface is
singular at $p_0$ and asymptotes a space-like cone with vertex at
$p_0$ and slope $\varphi_\infty$ where $\sinh\varphi_\infty =
2|H|/\sqrt{-\kappa}$.}

\subsection{Uniqueness of annular CMC surfaces}

We fix $\varepsilon=-1$ and $\kappa\le 0$ on this section. We then
present a version of a theorem proved by R. L\'opez (see \cite{Lo},
Theorem 1.2) about uniqueness of annular CMC in Minkowski space
$\mathbb{L}^3$.

Let $\Sigma_1$ be a  connected CMC  space-like surface in $\Mm\times
\Rr$ whose boundary is a geodesic circle $\Gamma$ in some plane
$\mathbb{P}_{a}$. We suppose that $\Sigma_1$ is a graph over
$\mathbb{P}_a-\Omega$, where $\Omega$ is the domain bounded by
$\Gamma$ on $\mathbb{P}_a$. We further suppose that the angle of
$\Sigma_1$ with respect to the planes $\Pp$ asymptotes, when
$\Sigma_1$ approaches $\partial_\infty \Mm\times\Rr$, a value
$\varphi_\infty^1$ so that $\sinh(\varphi_\infty^1)\ge
2|H|/\sqrt{-\kappa}$. We then consider $\Sigma_2$ a revolution
surface with same mean curvature, boundary and flux than $\Sigma_1$.
That this is possible we infer from the description on Theorem 5
above. From  the same theorem, we know that the asymptotic angle for
$\Sigma_2$ is
$\varphi_\infty^2=\textrm{arcsinh}(2|H|/\sqrt{-\kappa})$.


Suppose that $\Sigma_1\neq \Sigma_2$. Now, we move $\Sigma_1$
upwards until there is no contact with $\Sigma_2$. This is possible
since the asymptotic angle of $\Sigma_1$ is greater than or equal to
the asymptotic angle of $\Sigma_2$. Denote by $\Sigma_1 (t)$ the
copy of $\Sigma_1$ translated $t$ upwards (so that
$\Sigma_1(0)=\Sigma_1$). Then we define $t_0$ as the height where
occurs the first contact point. Suppose that $t_0>0$. Then, the
first contact is not at an interior point. Otherwise, by the
interior maximum principle, the surfaces are coincident, what
contradicts our hypothesis.
If the asymptotic angles are different, there are no point of
contact at infinity. If the angles are equal, then for small
$\delta$ the surfaces $\Sigma_1(t_0-\delta)$ and $\Sigma_2$
intersect transversally. We claim that there exists a connected
component $\Gamma'$ on $S=\Sigma_1(t_0-\delta)\cap\Sigma_2$ which is not null
homologous on both surfaces. Since both graphs have the topology of
a punctured plane, this means that $\Gamma'$ must be homologous to
$\Gamma$ on $\Sigma_2$.  Suppose by contradiction that all
components of $S$ will be null homologous. So, each  component
$\Gamma'$ of $S$ bounds a disc on both the graphs with common
boundary given by $\Gamma'$. These two discs are graphs over a disc
on $\mathbb{P}_a$ with the same mean curvature and same boundary. By
maximum principle they are equal. By analyticity, this implies that
the graphs coincide globally. From this contradiction, we conclude
that there exists component $\Gamma'$ of $S$ not null homologous.
The flux of $\Sigma_1(t_0-\delta)$ and $\Sigma_2$ through $\Gamma'$
are both equal to the flux of $\Sigma_1$ and $\Sigma_2$ through
$\Gamma$. However, after crossing $\Sigma_2$ along $\Gamma'$ towards
$\partial_\infty \Mm\times \Rr$, the surface $\Sigma_1(t_0-\delta)$
remains below $\Sigma_2$.
Then since $\partial_t$ is a time-like vector, it holds that
\[
\ae \eta_2, \partial_t\ad < \ae \eta_1, \partial_t\ad.
\]
along $\Gamma'$, where $\eta_1$ and $\eta_2$ are the {\it outward}
unit co-normal  of $\Sigma_a(t_0-\delta)$ and $\Sigma_2$ along
$\Gamma'$. However, this contradicts the  fact that the flux is the
same on both surfaces. This contradiction implies that $t_0=0$. Now,
if the surfaces contact at the boundary, they coincide globally, by
the boundary maximum principle. If not, then the angles satisfy
again a strict inequality and therefore the flux is not the same for
the two surfaces, a contradiction. We conclude from these
contradictions that $\Sigma_1=\Sigma_2$.


\vspace{0.3cm} \noindent {\bf Theorem 6. }{\it Let $\Sigma$ be a
space-like CMC surface on $\Mm\times \Rr$, $\kappa\le 0$, whose
boundary is a geodesic circle on a horizontal plane $\mathbb{P}_a$.
We suppose that $\Sigma$ is a  graph over the domain in
$\mathbb{P}_a$ outside the disc bounded by $\partial\Sigma$. We
further suppose that the angle between $\Sigma$ and the horizontal
planes asymptotes $\varphi_\infty$ with $\varphi_\infty\ge
\textrm{arcsinh} (2|H|/\sqrt{-\kappa})$. Then, $\Sigma$ is contained
on a revolution surface whose axis passes through the center of
$\partial\Sigma$ on $\mathbb{P}_a$.}

\vspace{0.3cm}

A similar reasoning shows, under the same hypothesis on the
asymptotic angle, that an entire space-like surface with an isolated
singularity and constant mean curvature  is a singular revolution
surface (v. \cite{Lo}, Theorem 1.3).

\section{Hopf differentials in some product spaces}

Let $\Sigma$ be a Riemann surface and $X:\Sigma\to \Mm\times \Rr$ be
an
isometric immersion. If $\kappa\ge 0$, we may consider $\Sigma$ as
immersed in $\Rr^4=\Rr^3\times \Rr$. If $\kappa <0$, we immerse
$\Sigma$ in $\mathbb{L}^{3}\times\Rr$. In fact, we may write
$X=(p,t)$, with $t\in \Rr$ and $p\in\Mm\subset\Rr^3$, in the first
case and $p\in \Mm \subset \mathbb{L}^{3}$ for $\kappa<0$. By
writing $\Mm\times\Rr\subset \mathbb{E}^4$ we mean all these
possibilities. The metric and covariant derivative in $\mathbb{E}^4$
are also denoted by $\langle\cdot,\cdot\rangle$ and $D$
respectively. We denote by $\epsilon$ the sign of $\kappa$. Recall
that $\varepsilon=1$ for Riemannian products and $\varepsilon=-1$
for Lorentzian ones.

Let $(u,v)$ be local coordinates in $\Sigma$ for which $X(u,v)$ is a
conformal immersion inducing the metric $e^{2\omega}\,(\dd u^2+ \dd
v^2)$ in $\Sigma$. So, denote by $\partial_u,
\partial_v$ the coordinate vectors and let $e_1=
e^{-\omega}\partial_u, \, e_2 = e^{-\omega}\partial_v$ be the
associated local orthonormal frame tangent to $\Sigma$. The unit
normal directions to $\Sigma$ in $\mathbb{E}^4$ are denoted by
$n_1,n_2=p/r$, where $r=(\epsilon\,\langle p,p\rangle)^{1/2}$. We
denote by $h^k_{ij}$ the components of $h^k$, the second fundamental
form of $\Sigma$ with respect to $n_k$, $k=1,2$. Then
\[
h^k_{ij}=\langle D_{e_i}e_j, n_k \rangle.
\]
It is clear that the $h^1_{ij}$ are the components of the second
fundamental form of the immersion $\Sigma\looparrowright \Mm\times
\Rr$. The components of $h^2$ are
\begin{eqnarray*}
& & h^2_{ij}= \langle D_{e_i}e_j, n_2\rangle= \langle
D_{e_i^h}e_j^{h},p/r\rangle= -\frac{1}{r}\langle
e_i^{h},e_j^{h}\rangle=
 \frac{1}{r}\big(\varepsilon\langle e_i^t,e_j^t\rangle- \delta_{ij}
\big)\\
& & =\frac{1}{r}\big(\varepsilon\langle e_i,\partial_t \rangle
\langle e_j,\partial_t
\rangle-\delta_{ij}\big)=\frac{\varepsilon}{r}\,\langle
e_i,\partial_t \rangle \langle e_j,\partial_t
\rangle-\frac{1}{r}\delta_{ij}.
\end{eqnarray*}
We remark that $\kappa = \epsilon/r^2$.  The components of $h^1$ and
$h^2$ in the frame $\partial_u,\partial_v$ are respectively
\begin{eqnarray*}
& &e= h^1(\partial_u, \partial_u)= e^{2\omega}h^1_{11},\, f
=h^1(\partial_u, \partial_v)= e^{2\omega}h^1_{12},\, g=
h^1(\partial_v, \partial_v)= e^{2\omega}h^1_{22}
\end{eqnarray*}
and
\begin{eqnarray*}
& & \tilde e= h^2(\partial_u, \partial_u)= e^{2\omega}h^2_{11},\,
\tilde f= h^2(\partial_u, \partial_v)= e^{2\omega}h^2_{12}, \,
\tilde g = h^2(\partial_v, \partial_v)= e^{2\omega}h^2_{22}.
\end{eqnarray*}
The Hopf differential associated to $h^k$ is defined by $\Psi^k=
\psi^k \dd z^2$, where $z = u+ iv$ and the coefficients $\psi^1,
\psi^2$ are
\[
\psi^1 = \frac{1}{2}(e-g)- i\, f,\quad \psi^2 = \frac{1}{2}(\tilde
e-\tilde g)- i\, \tilde f.
\]
The mean curvature of $X$ is by definition $H = (h^1_{11}+
h^1_{22})/2$. Differentiating the real part of $\psi^1$ we obtain
\begin{eqnarray*}
& & \partial_u \Big(\frac{e-g}{2}\Big)= \partial_u
\Big(\frac{e+g}{2}-g\Big)=\partial_u (e^{2\omega}H)-\partial_u
g=\partial_u (e^{2\omega}H)-\partial_u
\big(h^1(\partial_v,\partial_v)\big)\\
& & =\partial_u (e^{2\omega}H)-\big(D_{\partial_u}
h^1(\partial_v,\partial_v)+2h^1(D_{\partial_u}\partial_v,\partial_v)
\big)\\
& & =\partial_u (e^{2\omega}H) -
\big(D_{\partial_v}h^1(\partial_u,\partial_v)+\langle \bar
R(\partial_u,\partial_v)n_1,\partial_v\rangle+2h^1(D_{\partial_u}\partial_v,\partial_v)\big)\\
& & =\partial_u (e^{2\omega}H) - \big(\partial_v
(h^1(\partial_u,\partial_v))
-h^1(D_{\partial_v}\partial_u,\partial_v)-h^1(\partial_u,D_{\partial_v}\partial_v)\\
& &+\langle \bar
R(\partial_u,\partial_v)n_1,\partial_v\rangle+2h^1(D_{\partial_u}\partial_v,\partial_v)\big)\\
& &+\langle \bar
R(\partial_u,\partial_v)n_1,\partial_v\rangle\big)\\
& & =\partial_u (e^{2\omega}H) - \big(\partial_v f +\Gamma_{12}^1 f
+\Gamma_{12}^2 g-\Gamma_{22}^1 e-\Gamma^2_{22}f +\langle \bar
R(\partial_u,\partial_v)n_1,\partial_v\rangle\big)\\
& & =\partial_u (e^{2\omega}H) - \big(\partial_v f+f\partial_v\omega
+g\partial_u\omega +e\partial_u\omega -f\partial_v\omega +\langle
\bar
R(\partial_u,\partial_v)n_1,\partial_v\rangle\big)\\
& & =\partial_u (e^{2\omega}H) - \big(\partial_v f
+(e+g)\partial_u\omega  +\langle \bar
R(\partial_u,\partial_v)n_1,\partial_v\rangle\big)\\
& & = \partial_u (e^{2\omega}H)-2e^{2\omega}H\partial_u
\omega-\partial_v f -\langle \bar
R(\partial_u,\partial_v)n_1,\partial_v\rangle\\
& & =-\partial_v f +e^{2\omega}\partial_u H -\langle\bar
R(\partial_u,\partial_v)n_1,\partial_v\rangle.
\end{eqnarray*}
By similar calculations we also obtain
\[
\partial_v \Big(\frac{e-g}{2}\Big)= \partial_u f  - e^{2\omega}
\partial_v H +\ae \bar R(\partial_v, \partial_u)n_1, \partial_u \ad.
\]
We used above the Codazzi equation
\[
D_{\partial_u}h^1(\partial_v,\partial_v)=D_{\partial_v}h^1(\partial_u,\partial_v)+\langle
\bar R(\partial_u,\partial_v)n_1,\partial_v\rangle
\]
and the following expressions for the Christoffel symbols
$\Gamma^k_{ij}$ for the metric $e^{2\omega}\delta_{ij}$ in $\Sigma$
\[
\Gamma^1_{11}=-\Gamma^1_{22}=\Gamma^2_{12}=\partial_u \omega,\,
\Gamma^2_{22}=-\Gamma^2_{11}=\Gamma^1_{12}=\partial_v\omega.
\]
An easy calculation yields the components of the curvature tensor
\begin{eqnarray*}
\langle \bar R(\partial_u, \partial_v)n_1, \partial v\rangle =
\kappa\, e^{2\omega}\langle
\partial_u^h, n_1^h \rangle,\quad \langle \bar R(\partial_v, \partial_u)n_1, \partial u\rangle =
\kappa\,e^{2\omega}\langle
\partial_v^h, n_1^h\rangle.
\end{eqnarray*}
By this way, we then obtain the following pair of equations
\begin{eqnarray} & &\partial_u
\Re \psi^1 = \partial_v \Im \psi^1 -\kappa\, e^{2\omega}\langle
\partial_u^h, n_1^h \rangle + e^{2\omega} \partial_u H,\\
& & \partial_v \Re\psi^1 = -\partial_u \Im \psi^1 +\kappa\,
e^{2\omega}\langle
\partial_v^h, n_1^h \rangle - e^{2\omega} \partial_v H.
\end{eqnarray}
One also calculates
\begin{eqnarray*}
& & \partial_u \Re \psi^2 = \frac{\varepsilon}{2r}\, \partial_u
\Big(\langle
\partial_u,\partial_t\rangle^2
-\langle\partial_v,\partial_t\rangle^2\Big)= \frac{\varepsilon}{r}
\Big(\ae
\partial_u,\partial_t\rangle\,\ae D_{\partial_u}\partial_u,
\partial_t\ad \\
& & - \langle\partial_v,\partial_t\rangle\, \ae
D_{\partial_u}\partial_v, \partial_t\ad\Big)=
\frac{\varepsilon}{r}\Big(\ae
\partial_u,\partial_t\rangle\,\ae D_{\partial_u}\partial_u,
\partial_t\ad - \langle\partial_v,\partial_t\rangle\, \ae
D_{\partial_v}\partial_u, \partial_t\ad\Big)\\
& & =\frac{\varepsilon}{r} \Big(\ae
\partial_u,\partial_t\rangle\,\ae D_{\partial_u}\partial_u,
\partial_t\ad - \partial_v(\langle\partial_v,\partial_t\rangle\,  \ae
\partial_u, \partial_t\ad )+\ae D_{\partial_v}\partial_v, \partial_t\ad
\,
\ae \partial_u,\partial_t\ad \Big)\\
& & =\frac{\varepsilon}{r} \ae \partial_u,\partial_t\ad\,\ae
D_{\partial_u}\partial_u +D_{\partial_v}\partial_v ,
\partial_t\ad -
\frac{\varepsilon}{r}\,\partial_v\big(\ae\partial_u,\partial_t\ad\,\langle\partial_v,\partial_t\rangle\big)\\
& & = \frac{1}{r}\,\ae \partial_u, \partial_t \ad e^{2\omega}\,
\Delta t
-\frac{\varepsilon}{r}\,\partial_v\big(\ae\partial_u,\partial_t\ad\,\langle\partial_v,\partial_t\rangle\big)
\\
& & = 2H \frac{1}{r}e^{2\omega}\ae
\partial_u,\partial_t\ad\,\ae n_1 ,
\partial_t\ad -\frac{\varepsilon}{r}
\,\partial_v\big(\ae\partial_u,\partial_t\ad\,\langle\partial_v,\partial_t\rangle\big)=
-2H \frac{\varepsilon }{r}e^{2\omega}\ae
\partial_u^h,n_1^h\ad \\
& & -\frac{\varepsilon}{r}
\,\partial_v\big(\ae\partial_u,\partial_t\ad\,\langle\partial_v,\partial_t\rangle\big)
= -2H \frac{\varepsilon }{r}\, e^{2\omega} \ae
\partial_u^h,n_1^h\ad +
\,\partial_v\Im\psi^2.
\end{eqnarray*}
We used above the formula $\Delta t = 2H \ae n_1,
\partial_t\ad$, where $\Delta$ is the Laplacian on $\Sigma$ (see Section 6).
Similarly, we prove that
\[
\partial_v \Re \psi^2 = -\partial_u \Im \psi^2 +2H\frac{\varepsilon }{r}\, e^{2\omega}\ae
\partial_v^h,n_1^h\ad.
\]
Then, using the above mentioned fact that $\kappa = \epsilon/r^2$,
we conclude that the function $\psi:= 2H\psi^1 - \varepsilon
\frac{\epsilon}{r}\,\psi^2$ satisfies
\begin{eqnarray*}
& & \partial_u \Re \psi=\partial_v \Im  \psi + 2\Re \psi^1
H_u-2\Im\psi^1 H_v + 2e^{2\omega}HH_u
= \partial_v \Im  \psi + 2eH_u+2f H_v, \\
& & \partial_v \Re \psi=-\partial_u \Im \psi + 2\Re \psi_1
H_v+2\Im\psi^1 H_u -2e^{2\omega}HH_v = -\partial_u \Im \psi  - 2gH_v
-2f H_u.
\end{eqnarray*}
Now, using the complex parameter $z=u+iv$ and the complex derivation
$\partial_{\bar{z}} = \frac{1}{2}(\partial_u +
 i\partial_v) $  we get
\begin{eqnarray*}
\partial_{\bar{z}} \psi & = & (\partial_u \Re \psi - \partial_v \Im
\psi) +i (\partial_v \Re \psi + \partial_u \Im \psi)
 \\
& =  & 2eH_u+2fH_v -2if H_u-2igH_v
\end{eqnarray*}
%
That is,  defining the quadratic differential $Q:=2H\, \Psi^1
-\varepsilon\frac{\epsilon}{r} \, \Psi^2$ we prove that $Q$ is
holomorphic on $\Sigma$ if $H$ is constant. Inversely, if $Q$ is
holomorphic then
\[
e\, H_u +f\, H_v=0, \quad f\, H_u+  g\, H_v =0
\]
We may write this system in the following matrix form
\begin{eqnarray*}
\left[\begin{array}{ll} e & f \\
f & g \end{array}\right]\,\left[\begin{array}{ll} H_u\\
H_v\end{array}\right]=\left[\begin{array}{ll} 0\\
0\end{array}\right].
\end{eqnarray*}
This implies that $A\nabla H=0$, where $A=\ae \dd X, \dd
X\ad^{-1}\ae \dd n_1, \dd X\ad $ is the shape operator for $X$ and
$\nabla H$ is the gradient of $H$ on $\Sigma$. If $\nabla H=0$,
i.e., $H_u=H_v=0$ on $\Sigma$, then $H$ is constant. Thus, we may
suppose that $\nabla H\neq 0$ on an (open) set $\Sigma'$ of
$\Sigma$. On $\Sigma'$ we have $K_{\textrm{ext}}=:\det A=0$ .
However, $\det A=0$ is a closed condition. So, $\Sigma'$ is clopen
and therefore $\Sigma'=\Sigma$. Thus, $e_1=:\nabla H/|\nabla H|$ is
a principal direction with principal curvature $\kappa_1=0$.
Moreover $H=\kappa_2$, where $\kappa_2$ is the principal curvature
of $\Sigma$  calculated on a direction $e_2$ perpendicular to $e_1$.
So, the only  planar (umbilical) points on $\Sigma$ are the points
where $H$ vanishes. Moreover, the integral curves of $e_2$ are level
curves for $H=\kappa_2$ since they are orthogonal to $\nabla H$.
Thus, $H$ is constant along such each line. So,  we  proved

\vspace{0.3cm}

\noindent{\bf Theorem 7.} {\it   The quadratic differential $Q=2H
\Psi^1-\varepsilon\frac{\epsilon}{r}\, \Psi^2$ is holomorphic on
$\Sigma$ if $H$ is constant. Inversely, if we suppose that $\Sigma$
is compact (more generally, if $\Sigma$ does not admit a function
without critical points, or a vector field without singularities),
then $H$ is constant if $Q$ is holomorphic.}

\vspace{0.3cm}

The considerations above imply that if there exist examples of
surfaces with holomorphic $Q$ and non constant mean curvature, these
examples must be non compact, have zero extrinsic Gaussian curvature
and are foliated by curvature lines along which $H$ is constant.
Recently, P. Mira and I. Fern\'andez announced to the authors had
constructed such examples.

For  $\varepsilon=1$, the quadratic form $Q$ coincides with that one
obtained by U. Abresch and H. Rosenberg in (\cite{AR}). It is clear
that $Q$ is the complexification of the traceless part of the second
fundamental form $q$ corresponding to the normal direction $2H n_1
-\varepsilon \frac{\epsilon}{r}\,n_2$ on the normal bundle of
$\Sigma \looparrowright \mathbb{E}^4$.

Using the Theorem 7, we present the following generalization of the
theorem of Abresch and Rosenberg quoted in the Introduction:

\vspace{0.3cm}

\noindent{\bf Theorem 8.} {\it Let $X:\Sigma\to \Mm\times\Rr$ be a
complete CMC immersion of a surface $\Sigma$ in $\Mm\times \Rr$. If
$\varepsilon=1$ and $\Sigma$ is homeomorphic to a sphere, then
$X(\Sigma)$ is a rotationally invariant spherical surface. If
$\Sigma$ is homeomorphic to a disc and $Q\equiv 0$ on $\Sigma$, then
$X(\Sigma)$ is a rotationally invariant disc. For $\varepsilon=-1$
and $\kappa\le 0$, if $X(\Sigma)$ is simply-connected, space-like
and $Q\equiv 0$ on $\Sigma$, then the same conclusion holds.}

\vspace{0.3cm}

\noindent{\it Proof of the Theorem 8.} By hypothesis, we have
$Q\equiv 0$ (if $\Sigma$ is homeomorphic with a sphere, this follows
from the fact that $Q$ is holomorphic). Thus, $2H\psi^1 \equiv
\varepsilon\frac{\epsilon}{r}\,\psi^2$. Given an arbitrary local
orthonormal frame field $\{e_1,e_2\}$, we may write this as
\begin{eqnarray}
\label{eq1} & & 2H h^1_{12}= \kappa\,\langle e_1,\partial_t \rangle
\langle
e_2,\partial_t\rangle,\\
\label{eq2} & & 2H(h^1_{11}-h_{22}^1)= \kappa\,\langle
e_1,\partial_t\rangle^2 - \kappa\,\langle e_2,\partial_t\rangle^2.
\end{eqnarray}
If $H=0$, then it follows from these equations that the vector field
$\partial_t$ is always normal to $\Sigma$. So, the surface is part
of a plane $\Pp=\Mm\times\{t\}$, for some $t\in \Rr$. Since $\Sigma$
is complete, we conclude that $\Sigma=\Pp$.

We then may consider only  CMC surfaces with $H\neq 0$. If $(p,t)$
is an  umbilical point of $\Sigma$ we have  for an arbitrary frame
that $h^{1}_{12}=0$ at this point. So, either $\ae e_1,
\partial_t\ad =0$ or $\ae e_2,\partial_t\ad=0$ at $(p,t)$. Since
 $h^1_{11}=h^1_{22}=H$ at $(p,t)$ the equation (\ref{eq2})
implies that both angles $\ae e_i, \partial_t\ad$ are null. So, we
conclude that if $Q=0$, then umbilical points are the points where
$\Sigma$ has horizontal tangent plane, and vice-versa.

If $(p,t)$ is not an umbilical point in $\Sigma$, we may choose the
frame $\{e_1,e_2\}$ as principal frame locally defined (on a
neighborhood $\Sigma'$ of that point). Thus, $h^1_{12}=0$ and
therefore $\langle e_1,\partial_t\rangle =0$ or $\langle e_2,
\partial_t\rangle =0$ on $\Sigma'$. We fix $\langle
e_1,\partial_t\rangle =0$. If we denote by $\tau$ the tangential
part $\partial_t-\varepsilon\ae \partial_t, n_1\ad n_1$ of the field
$\partial_t$, then $\tau= \ae e_2,
\partial_t\ad \, e_2$. Thus from (\ref{eq2})
it follows that the principal curvatures of $\Sigma$ are
\[
h^1_{11}=H-\frac{\kappa}{4H}\,|\tau|^2,\quad
h^1_{22}=H+\frac{\kappa}{4H}\, |\tau|^2.
\]
The lines of  curvature on $\Sigma'$ with direction $e_1$ are
locally contained in the planes $\Pp$. Inversely, the connected
components of $\Sigma' \cap \Pp$ are lines of curvature with tangent
direction given by $e_1$. Thus, if we parameterize such a line by
its arc length $s$, we have
\begin{eqnarray}
\frac{\dd}{\dd s} \langle n_1,\partial_t\rangle  = \langle D_{e_1}
n_1,
\partial_t\rangle = h^1_{11}\langle e_1, \partial_t\rangle = 0.
\end{eqnarray}
We conclude that, for a fixed $t$,  $\Sigma'$ and $\Pp$ make a
constant angle $\theta(t)$ along each connected component of their
intersection. So, if a connected component of the intersection
between $\Pp$ and $\Sigma$ has a non umbilical point, then the angle
is constant, non zero, along this component, unless that there
exists also an umbilical point on this same component. However at
this point the angle is necessarily zero. So, by continuity of the
angle function, either all points on a connected component
$\Sigma\cap \Pp$ are umbilical and the angle is zero, or all points
are non umbilical and the angle is non zero. However, supposes that
all points on a connected component $\sigma$ are umbilical points
for $h^1$. Then, as we noticed above, $\Sigma$ is tangent to $\Pp$
along $\sigma$. So, along $\sigma$, we have $\ae e_1, \partial_t\ad
= \ae e_2,
\partial_t\ad=0$ and therefore by equations (\ref{eq1}) and (\ref{eq2}) we have
$h^1_{ii}=0$ and $H=0$. From this contradiction, we conclude that
the umbilical points may not be on {\it any} curve on $\Sigma \cap
\Pp$. The only possibility is that there exist isolated umbilical
points as may occurs on the top and bottom levels $t=a$ and $t=b$ of
$X(\Sigma)$.

So, there exists an orthonormal principal frame field $\{e_1,e_2\}$
on a dense subset of $\Sigma$. On this dense subset we have
$\tau\neq 0$ and then we may choose a positive sign for $\sin
\theta(t)$ or $\sinh\,\theta(t)$, where $\theta(t)$ is the angle
between $n_1$ and $\partial_t$ along a given component of
$\Sigma\cap \Pp$. We denote both of these functions by the same
symbol $\sn(t)$. Now, we calculate the geodesic curvature of the
horizontal curvature lines on $\Pp$. We have
\[
e_2 = \frac{\tau}{|\tau|}= \frac{1}{\sn(t)}\,\tau=
\frac{1}{\sn(t)}\,(\partial_t -\varepsilon\ae \partial_t, n_1\ad
n_1)=\frac{1}{\sn(t)}\,(\partial_t -\dot\sn(t) n_1)
\]
Since $\ae n_1,
\partial_t \ad$ is constant along this curve and therefore $\sn(t)$
is constant we conclude that
\[
D_{e_1}e_2 = \frac{1}{\sn(t)}\big(D_{e_1}\partial_t - \dot\sn(t) \,
D_{e_1}n_1\big)=\frac{\dot\sn(t)}{\sn(t)}\, h^1_{11}e_1
\]
where $\dot \sn(t)=\cos\theta(t)$ for $\varepsilon=1$ and $\dot
\sn(t)=\cosh\theta(t)$ for $\varepsilon=-1$. So the geodesic
curvature $\ae D_{e_1}e_1, e_2\ad$ of the horizontal lines of
curvature relatively to $\Sigma$  is given by
$-(\dot\sn(t)/\sn(t))\, h^1_{11}$. This means that the horizontal
lines of curvature  have constant geodesic curvature on $\Sigma$.
Now, defining $\nu = Je_1  = \varepsilon\sn(t)\, n_1 -\dot \sn(t)\,
e_2$, we calculate
\[
\ae D_{e_1}\nu, e_1 \ad = -\varepsilon\sn(t) h^1_{11}
-\dot\sn(t)\frac{\dot\sn(t)}{\sn(t)}\, h^1_{11}= -\frac{1}{\sn(t)}\,
h^1_{11}.
\]
Thus,  it follows that the geodesic curvature of the horizontal
lines of curvature on $\Sigma \cap \Pp$ relatively to the plane
$\Pp$ is also constant and equal to $h^1_{11}/\sn(t)$. We conclude
that for each $t$, $\Sigma\cap \Pp$ consists of constant geodesic
curvature lines of $\Pp$.

We also obtain $\ae D_{e_2}e_2,e_1\ad = 0$. So, the curvature lines
of $\Sigma$ with direction $e_2$ are geodesics on  $\Sigma$. We then
prove that these lines are contained on vertical planes.  Fixed a
point $(p,t)$ in $\Sigma\cap \Pp$, let $\alpha(s)$ be the line of
curvature with $\alpha'=e_2$ passing by $(p,t)$ at $s=0$. We want to
show that $\alpha$ is contained on the vertical geodesic plane $\Pi$
determined by $e_2(p,t)$ and $\partial_t$. This is the plane spanned
by $e_2$ and $n_1$ at $(p,t)$. For each $s$, consider the vertical
geodesic plane $\Pi_s$ on $\Mm\times \Rr$ for which $e_2=\alpha'(s)$
and $D_{e_2}e_2=D_{\alpha'}\alpha'$ are tangent at $\alpha(s)$. This
plane is of the form $\sigma_s\times \Rr$, where $\sigma_s$ is some
geodesic on $\Mm$ which by its turn is the intersection of $\Mm$ and
some plane $\pi_s$ on $\mathbb{E}^3$ with unit normal $a(s)$. The
intersection of the hyperplane $\pi_s \times \Rr$ of $\mathbb{E}^4$
with $\Mm\times \Rr$ is then the plane $\Pi_s$. Now $p(s)\wedge
\alpha'(s)\wedge D_{\alpha'}\alpha'$ is a normal direction  to that
hyperplane on $\mathbb{E}^4$ where $p(s)=\alpha(s)^h$. However,
since $\alpha$ is at the same time line of curvature and geodesic
then
\[
D_{\alpha'}\alpha' = D_{e_2}e_2 = (D_{e_2}e_2)^T +(D_{e_2}e_2)^N =
(D_{e_2}e_2)^N = h_{22}^1\, n_1.
\]
Thus we conclude that the unit normal to the hyperplane $\Pi_s$ is
\[
a(s)= p(s) \wedge e_2(s)\wedge n_1(s).
\]
Differentiating  we obtain $a'=0$.
So $a(s)$ is constant. Thus implies that  $\Pi_s=\Pi$ for all $s$.
So, $\alpha(s)$ is a plane curve contained in $\Pi$. Notice that
$\Pi$ has normal $e_1(p,t)$ since $e_1(p,t)=a(0)$. We then conclude
that the integral curves of $e_2$ are planar geodesics on $\Sigma$.

So, for a fixed $t$, let $\sigma(s)$ be a component of $\Sigma \cap
\Pp$. Then $\sigma$ is a constant geodesic curvature curve on $\Pp$.
Moreover,  the vertical plane passing through $\sigma(s)$ with
normal $e_1(\sigma(s))$ is a symmetry plane of $\Sigma$ since
contains a geodesic of $\Sigma$, namely the curvature line in
direction $e_2$ passing through $\sigma(s)$. Thus, the surface is
invariant with respect to the isometries fixing $\sigma$. Since the
surface is homeomorphic to a disc or a sphere (see Remark 2 below),
then we conclude that these isometries are elliptic (their orbits
are closed circles). This means that $X(\Sigma)$ is rotationally
invariant in the sense of Section 1. So, the proof is concluded.

\vspace{0.3cm}

\noindent {\it Remark 1.} We also prove the Theorem 8 by the
following reasoning: denote by $\Pi_s$ the plane passing through
$\sigma(s)$ with normal $e_1$. This plane contains the curvature
line with initial data $\sigma(s)$ for position and $e_2(\sigma(s))$
for velocity. Its plane curvature is given by the derivative of its
angle with respect to the (fixed) direction $\partial_t$, that is,
$\theta(t)$. These data, by the fundamental theorem on planar
curves, determine completely the curve. Changing the point on
$\sigma$, the initial data differ by a rigid motion  (an isometry on
$\Pp$) and the curvature function remains the same at points of
equal height. Then, by the uniqueness part on the theorem cited
before, the two curves differ only by the same rigid motion. This
means that the surface is invariant by the rigid motions fixing
$\sigma$. Thus, the proof is finished by proving that the only
possible isometries are the elliptic ones.



\noindent{\it Remark 2.} For $\kappa\le 0$ and $\varepsilon=-1$,
since $X(\Sigma)$ is space-like, it is acausal. Thus, the coordinate
$t$ is bounded on $\Sigma$. Moreover, the projection $(p,t)\in
\Sigma\mapsto p\in \Mm$ increases Riemannian distances. So, is a
covering map and therefore $X(\Sigma)$ is locally a graph over the
horizontal planes. If we suppose $\Sigma$ simply connected, then
$X(\Sigma)$ is globally diffeomorphic with  $\Pp$. Is, in fact, a
disc-type graph.

\vspace{0.5cm}

Let $X:\Sigma\to \Mm\times \Rr$ be an  immersion  of a surface with
boundary. We suppose that $X|_{\partial\Sigma}$ is a diffeomorphism
onto its image $\Gamma=X(\partial\Sigma)$. We further suppose that
$X(\partial\Sigma)$ is contained on some plane $\Pp$. So, $\Gamma$
is a embedded curve  on $\Pp$ that bounds a domain $\Omega$. In what
follows we always make this hypothesis while treating immersions of
surfaces with boundary. Now, we fix $\varepsilon=-1$ and suppose
that $X(\Sigma)$ is space-like. We may prove under these assumptions
that $\Sigma$ is simply-connected (disc-type) and $X(\Sigma)$ is a
graph over $\Omega$. This conclusion also holds if $\Gamma$ is
supposed to be a graph over some embedded curve on $\Pp$.

Thus, if we suppose either $\varepsilon=1$ and $\Sigma$ a disc, or
$\varepsilon=-1$  (with the additional hypothesis that $Q=0$ on both
cases) then we are able to prove that if $X(\Sigma)$ is an immersed
CMC surface with boundary, then $X(\Sigma)$ is contained on a
rotationally invariant CMC disc. In fact, the reasoning on Theorem 8
works well on these cases to show that $X(\Sigma)$ is foliated by
geodesic circles and that the angle with a plane $\Pp$ is constant
along $\Sigma \cap \Pp$. This suffices to show that $X(\Sigma)$ is
rotationally invariant.

\section{Free boundary surfaces in product spaces}

A classical result  of J. Nitsche (see, e.g., \cite{Ni}, \cite{RV}
and \cite{So}) characterizes discs and spherical caps as equilibria
solutions for the free boundary problem in space forms. We will be
concerned now about to reformulate this problem in the product
spaces $\Mm\times \Rr$.

Let $\Sigma$ be an orientable compact surface with non empty
boundary and $X: \Sigma \to \Mm\times \Rr$ be an isometric
immersion. By a \emph{volume-preserving variation of $X$} we mean a
family $X_s: \Sigma\to \Mm\times \Rr$ of isometric immersions such
that $X_0 = X$ and $\int \langle
\partial_s X_s, n_s\rangle \dd A_s=0$, where $\dd A_s$ and $n_s$
represent respectively the element of area  and an unit normal
vector field to $X_s$. In the sequel  we set $\xi=
\partial_s X_s$ and $f=\ae \xi_s,n_s\ad$ at $s=0$. We say that $X_s$ is an \emph{admissible
variation} if it is volume-preserving and at each time $s$ the
boundary $X_s(\partial \Sigma)$ of $X_s(\Sigma)$ lies on a
horizontal plane $\mathbb{P}_a$. We denote by $\Omega_s$ the compact
domain in $\mathbb{P}_a$ whose boundary is $X_s(\partial \Sigma)$
(in the spherical case $\kappa>0$, we choose one of the two domains
bounded by $X_s(\partial\Sigma)$). A {\it stationary surface} is by
definition a critical point for the following functional
\[
E(s) = \int_\Sigma \dd A_s + \alpha\int_{\Omega_s} \dd \Omega,
\]
for some constant $\alpha$, where $\dd \Omega$ is the volume element
for $\Omega_s$ induced from $\mathbb{P}_a$. The first variation
formula for this functional is (see \cite{RV} and \cite{AP} for the
corresponding formulae in space forms)
\[
E'(0)= -2\int_\Sigma H f + \int_{\partial \Sigma}\ae \xi,
\eta+\alpha  \bar\eta\ad \,\dd \sigma,
\]
where $\dd \sigma$ is the line element for $\partial \Sigma$ and
$\eta,\, \bar \eta$ are the unit co-normal vector fields to
$\partial \Sigma$ relatively to $\Sigma$ and to $\mathbb{P}_a$. If
we prescribe $\alpha= -\cos\theta$ in the Riemannian case and
$\alpha=-\cosh \theta$ in the Lorentzian case, then we conclude that
a stationary surface $\Sigma$ has constant mean curvature and  makes
constant angle $\theta$ along $\partial\Sigma$ with the horizontal
plane.

In what follows, {\it spherical cap} means that the surface is a
part of a CMC revolution sphere bounded by some circle contained in
a horizontal plane and centered at the rotation axis. Similarly, the
term {\it hyperbolic cap} means a part of a CMC rotationally
invariant disc bounded by a horizontal circle centered at the
rotation axis. Granted this, we state the following theorem.

\vspace{0.3cm}

\noindent {\bf Theorem 9. } {\it Let $\Sigma$ be a  surface with
boundary and let $X: \Sigma\to \Mm\times \Rr$ be a {\rm stationary}
immersion for free boundary admissible variations whose boundary
lies in some plane $\mathbb{P}_a$.  If $\varepsilon =1$ and $\Sigma$
is disc-type, then $X(\Sigma)$ is a spherical cap. If $\varepsilon
=-1$, then $X(\Sigma)$ is a hyperbolic cap.}

\vspace{0.3cm}

The proof of Theorem 9 follows closely the guidelines of the proof
of the Nitsche's Theorem in $\Rr^3$ as we may found in \cite{Ni} and
\cite{RV}. Let $\Sigma$ denote the disc $|z|<1$ in $\Rr^2$, where
$z=u+iv$. If we put $\partial_z = \frac{1}{2}(\partial_u -
i\partial_v)$, then the $\mathbb{C}$-bilinear complexification of
$q$ satisfies
\[
q_{\mathbb{C}}(\partial_z, \partial_z)= q(\partial_u,
\partial_u) -q(\partial_v, \partial_v) -2 i q(\partial_u,
\partial_v) = 2Q(\partial_z, \partial_z).
\]
%
%
%
%
%
%
%
Now, since $X(\partial \Sigma)$ is contained in $\mathbb{P}_a$ then
$q(\tau, \eta)=0$ on $\partial \Sigma$. Here $\tau = e^{-\omega}
(-v\partial_u + u
\partial_v)$ is the unit tangent vector to $\partial\Sigma$ and
$\eta
=e^{-\omega} (u \partial_u + v
\partial_v)$ is the unit outward co-normal to $\partial \Sigma$. In
fact $h^2(\tau, \eta) = 0$ since that $\tau$ is a horizontal vector
and $h^1(\tau, \eta)=0$ since that $\partial\Sigma$ is a line of
curvature for $\Sigma$ by Joachimstahl's Theorem.

On the other hand, we have on $\partial\Sigma$ that
\[
0=q(\tau, \eta) =(u^2-v^2)\,q(\partial_u, \partial_v) - uv
\,q(\partial_u,\partial_u) +uv \,q(\partial_v, \partial_v)=\Im
\big(z^2 Q (\partial_z, \partial_z)\big)
\]
 From this we conclude that $\Im (z^2
Q)\equiv 0$ on $\partial\Sigma$. Since $z^2 Q$ is holomorphic on
$\Sigma$, then $\Im z^2 Q$  is harmonic. So, $\Im z^2 Q=0$ on
$\Sigma$ and therefore $z^2 Q \equiv 0$ on $\Sigma$. Hence, $Q\equiv
0$ on $\Sigma$. This implies that $X(\Sigma)$ is part of a CMC
revolution sphere or a CMC rotationally invariant disc. This
finishes the proof of the Theorem 9.

\vspace{0.4cm}

We obtain also a result about stable CMC discs in $\Mm\times \Rr$,
following ideas presented in \cite{ALP}. Here, stability for a CMC
surface $\Sigma$ means that the quadratic form
\[
J[f] = \varepsilon\,\int_\Sigma \Big(\Delta f + \varepsilon
\big(|A|^2 + \textrm{Ric}(n_1, n_1)\big)f\Big) \, f \, \dd A,
\]
is non-negative with respect to the all variational fields $f$
generating preserving-volume variations (see \cite{BC} and \cite{BO}
for the case $\kappa=0$). In the formula above, $\textrm{Ric}$ means
the Ricci curvature tensor of $\Mm\times \Rr$.

\vspace{0.3cm}

\noindent{\bf Theorem 10.} {\it Let $\Sigma$ be an immersed surface
with boundary and constant mean curvature $H$ in $\Mm\times\Rr$.
Suppose that $\partial \Sigma$ is  a geodesic circle in some plane
$\mathbb{P}_a$ and that the immersion is stable. For $\varepsilon=1$
we further suppose that $\Sigma$ is disc-type and for
$\varepsilon=-1$ that the immersion is space-like. Then $\Sigma$ is
a spherical or hyperbolic cap, if $H\neq 0$. If $H=0$ then $\Sigma$
is a totally geodesic disc.}

\vspace{0.3cm}

We consider the vector field $Y(t,p) = a\wedge
\partial_t \wedge p$, where $a$ is the vector in $\mathbb{E}^3$ perpendicular to the plane
where $\partial \Sigma$ lies. This is a Killing field in $\Mm\times
\Rr$. Then $f=\langle Y, n_1\rangle$ satisfies trivially $J[f]=0$.
Let $\eta$ be the exterior unit co-normal direction to $\Sigma$
along the boundary $\partial \Sigma$.

The normal derivative of $f$ along $\partial \Sigma$ is calculated
as
\begin{eqnarray*}
& &\eta(f)  = \eta\langle Y, n_1\rangle = \langle a\wedge
\partial_t\wedge D_{\eta} p, n_1\rangle +
\langle a\wedge
\partial_t\wedge p, D_{\eta} n_1\rangle\\
& & =  \langle a\wedge
\partial_t\wedge \eta, n_1\rangle + \langle a\wedge
\partial_t\wedge p, D_{\eta} n_1\rangle =-\langle a\wedge
\partial_t\wedge n_1, \eta\rangle + \langle a\wedge
\partial_t\wedge p, D_{\eta} n_1\rangle\\
& & =  \ae \tau, \eta\ad +\langle \tau, D_{\eta}n_1\rangle = \langle
\tau, D_{\eta}n_1\rangle = -h^1(\tau, \eta),
\end{eqnarray*}
where $\tau = a\wedge \partial_t\wedge p $ (the restriction of $Y$
to the boundary of $\Sigma$) is the tangent positively oriented unit
vector to $\partial \Sigma$. Since that $\langle \tau,
\partial_t\rangle =0$ and $\langle \tau, \eta\rangle=0$ it follows that
\[
h^2(\tau, \eta) = -\frac{1}{r}\,\langle\tau^h, \eta^h\rangle =0.
\]
This yields
\[
2H \, \eta(f) = -2H \, h^1(\tau,\eta) = - q(\tau, \eta).
\]
However, if $u,v$ denote the usual cartesian coordinates on $\Sigma$
then
\[
q(\tau, \eta) = e^{-2\omega}\, q(u\partial_u +v\partial_v,
-v\partial_u + u\partial_v) = -\Im (z^2 Q)
\]
on $\partial \Sigma$. We conclude that $2H\, \eta(f) = \Im (z^2 Q)$.
Proceeding as in (\cite{ALP}) we verify that $\eta(f)$ vanishes at
least three times. Applying   Courant's theorem on nodal domains
allows us to conclude that $f$ vanishes on the whole disc. So,
$X(\Sigma)$ is foliated by the flux lines of $Y$, i.e. by horizontal
geodesic circles centered at the same vertical axis. So, $X(\Sigma)$
is a spherical or hyperbolic cap as we claimed. This proves Theorem
10.

\section{Flux formula and Killing graphs}

\subsection{Flux formula}

Let $\Sigma$ be an immersed surface in $\Mm\times\Rr$ with constant
mean curvature $H$. We denote by $\textrm{Div}$ and $\textrm{div}$
respectively the divergence operator on $\Mm\times \Rr$ and on
$\Sigma$. Consider a Killing vector field $Y$ on $\Mm\times \Rr$.
Thus restricting $Y$ to $\Sigma$ one finds
\[
\textrm {div} Y =: \sum_{i}\ae D_{e_i}Y, e_i \ad =0.
\]
However using the decomposition $Y = Y^T + Y^N = Y^T +\varepsilon
\ae Y, n\ad \, n$ we obtain
\begin{eqnarray*}
& & \textrm{div} Y = \textrm{div} Y^T + \textrm{div} Y^N =
\textrm{div} Y^T + \varepsilon \ae Y,n \ad  \ae D_{e_i} n, e_i\ad =
\ae \nabla_{e_i} Y^T, e_i \ad\\
& &  -2H\varepsilon \ae Y,n \ad = \textrm{div} Y^T -2H \varepsilon
\ae Y, n\ad.
\end{eqnarray*}
Then by Stokes's Theorem on $\Sigma$
\[
0=\int_\Sigma \textrm{div} Y\, \dd A = \int_{\partial\Sigma }\ae Y,
\nu\ad\, \dd \sigma -2H \varepsilon \int_{\Sigma} \ae Y, n\ad \, \dd
A,
\]
where $\nu$ is the {\it outward} unit co-normal vector field along
$\partial\Sigma$. By this way we obtain the {\it first Minkowski
formula}
\begin{equation}
\int_\Sigma \,\varepsilon H  \ae Y, n\ad \, \dd A= \frac{1}{2}\,
\int_{\partial\Sigma} \ae Y, \nu \ad \, \dd \sigma.
\end{equation}
In the case where $\Sigma$ and $\Omega$ are homologous oriented
cycles on $\Mm\times \Rr$, we conclude from the formula
$\textrm{Div} Y =0$ and divergence theorem that
\[
\int_\Sigma \ae Y, n\ad \, \dd A + \int_\Omega \ae Y, n_\Omega \ad
\, \dd\Omega=0.
\]
We then obtain the {\it flux formula} for Killing vector fields:
\begin{equation}
\label{Killing} \int_{\partial\Sigma }\ae Y, \nu\ad\, \dd \sigma +
2H \varepsilon \,\int_\Omega \ae Y, n_\Omega \ad \, \dd\Omega =0.
\end{equation}

\subsection{Killing graphs and  height estimates}

Let $n$ be an unit normal vector field to $\Sigma\looparrowright
\Mm\times\Rr$. We next consider the function $\langle Y,n\rangle$.
Let $e_1,e_2$ be an adapted orthonormal moving frame with $\nabla
e_i=0$ at a point $(p,t)\in \Sigma$. We may suppose that $e_i$ is
principal at that point. We have for $v\in T\Sigma$ that
\begin{eqnarray*}
& & v \langle Y, n\rangle = \langle D_{v} Y, n\rangle + \langle Y,
D_v n\rangle = \langle D_v Y, n\rangle + \varepsilon\langle D_v
n,n\rangle
\langle n,Y\rangle +\langle (D_v n)^T,Y\rangle \\
& & =-\langle v, D_n Y\rangle +\langle (D_v n)^T,Y\rangle=-\langle
v, D_n Y+A(Y^T)\rangle.
\end{eqnarray*}
Hence
\[
\nabla \langle Y,n\rangle = -A(Y^T) - (D_n Y)^T=\big((D_{Y^T}n)-(D_n
Y)\big)^T
\]
The restriction of a Killing field to a CMC  surface is a Jacobi
field for it. Then we have
\begin{eqnarray*}
\big(\Delta +\varepsilon |A|^2 +
\varepsilon\textrm{Ric}(n,n)\big)\langle Y,n\rangle=0.
\end{eqnarray*}
We also compute \begin{eqnarray*} \textrm{Ric}(n,n)=\kappa
\varepsilon(1-\langle n,\partial_t\rangle^2).
\end{eqnarray*}

We suppose that the distribution spanned by the vectors orthogonal
to $Y$ is integrable (this is a weaker condition than to assume $Y$
is closed). Let $\mathbb{N}$ be the domain in $\Mm\times \Rr$ free
of singularities of $Y$. So, $\mathbb{N}$ is  foliated by surfaces
orthogonal to the flow lines of $Y$. Let $s$ be the flow parameter
on the flow lines of $Y$, so that each leaf is a level surface for
$s$. Taking $s$ as a coordinate on $\mathbb{N}$, it is clear that
$\partial_s=Y$. We also have
\[
\bar \nabla s =g^{ss}\partial_s = \frac{Y}{|Y|^2}:= f Y.
\]
Then the gradient of $s$ restricted to a surface  $\Sigma$ on
$\mathbb{N}$ is $\nabla s = f Y^T$  and its Laplacian is calculated
as
\begin{eqnarray*}
& & \Delta s = \langle \nabla_{e_i} f Y^T, e_i \rangle =\langle
D_{e_i} f Y^T, e_i \rangle = \langle \nabla f, e_i\rangle \langle
Y^T, e_i\rangle + f\langle D_{e_i}  Y^T, e_i \rangle\\
& & =\langle \nabla f,  Y^T\rangle + f\langle D_{e_i}  Y, e_i
\rangle -\varepsilon f\langle D_{e_i}  \langle Y,n\rangle n,e_i
\rangle=\langle \nabla f,  Y^T\rangle  -\varepsilon f\langle
Y,n\rangle \langle D_{e_i} n,e_i \rangle.
\end{eqnarray*}
Thus $\Delta s =2H\varepsilon f \langle Y,n\rangle +\langle \nabla
f, Y\rangle$. However, we easily see that the Killing equation
implies that the norm of $Y$ is conserved along the flow lines of
$Y$. Then $\ae \nabla |Y|, Y\ad=0$ and therefore $\ae \nabla f, Y\ad
=0$. So
\[
\Delta s =2H\varepsilon f \langle Y,n\rangle.
\]
We also have from Jacobi's equation
\[
\Delta \ae Y, n\ad = -\varepsilon\big(|A|^2+\textrm{Ric}(n,n)
\big)\, \ae Y,n\ad =-\varepsilon\big(|A|^2+\kappa\varepsilon\,(1-\ae
n,\partial_t\ad^2) \big)\, \ae Y,n\ad.
\]
We then fix $\varepsilon=1$. Suppose that $\Sigma$ has boundary on
the leaf $\Pi$ given by $s=0$ and that $\ae Y,n\ad\ge 0$ on
$\Sigma$. So, $H\le 0$ when we consider $n$ pointing outwards $\Pi$.
Next, for a given constant $c$, the function $\phi=: Hc\, s +\ae
Y,n\ad$ satisfies $\phi|_{\partial\Sigma}\ge 0$ and
\begin{eqnarray*}
& & \Delta \phi= \big(2H^2 \frac{c}{|Y|^2}-|A|^2
-\kappa(1-\nu^2)\big)\, \ae Y,n\ad,
\end{eqnarray*}
where $\nu=\ae n,\partial_t\ad$. We want to choose $c$ so that
$\phi$ is super-harmonic. It suffices that
\begin{equation}
\label{c} 2H^2 \frac{c}{|Y|^2}-|A|^2 -\kappa(1-\nu^2)\le 0.
\end{equation}
However
\begin{eqnarray*}
& & |A|^2=k_1^2+k_2^2  = (k_1+k_2)^2 +(k_1-k_2)^2 -k_1^2-k_2^2=4H^2
+4|\psi^1|^2-|A|^2.
\end{eqnarray*}
So, $|A|^2 =2H^2 +2|\psi^1|^2$.
For further reference we point that (for any sign of $\varepsilon$)
\[
|\psi^2|^2 = \frac{\kappa^2}{4}\, (1-\nu^2)^2.
\]
Thus, (\ref{c}) is rewritten as
\[
2H^2 \big(\frac{c}{|Y|^2}-1\big)\le 2|\psi^1|^2 +\kappa(1-\nu^2)
\]
If $\kappa\ge 0$ it suffices to take $0<c\le \inf_\Sigma|Y|^2$. For
$\kappa<0$,
if we rather suppose that $2H^2+\kappa>0$, then we obtain the
super-harmonicity of $\phi$ when
\[
0<c\le \inf_\Sigma |Y|^2\frac{H^2+\frac{\kappa}{2}}{H^2}.
\]
Thus, for these choices for $c$ we have $\Delta\phi\le 0$ and
\begin{eqnarray}
\label{sestimates} s\le\frac{1}{|H|}\, \frac{\sup_\Sigma
|Y|}{\inf_\Sigma |Y|^2}, \, (\kappa>0),\quad  s\le
\frac{|H|}{H^2+\frac{\kappa}{2}}\,
\frac{\sup_\Sigma|Y|}{\inf_\Sigma|Y|^2},\, (\kappa<0).
\end{eqnarray}

\vspace{0.3cm} \noindent {\bf Theorem 11. }{\it Let $Y$ be a Killing
field on $\Mm\times \Rr$, $\varepsilon=1$, which determines an
integrable orthogonal distribution $\mathcal{D}$. Let $\Sigma$ be an
immersed CMC surface on $\mathbb{N}$ whose boundary lies on a
integral leaf  of $\mathcal{D}$. If $s$ is the parameter of the flow
lines of $Y$, then it holds the estimates on (\ref{sestimates}). If
$\Sigma$ is a compact closed embedded CMC surface on $\mathbb{N}$,
then $\Sigma$ is symmetric with respect to some integral leaf of
$\mathcal{D}$.}

\vspace{0.3cm}

The proof of the second statement on the theorem above is similar to
that one presented in Proposition 1 of \cite{HLR}. It is based on
Aleksandrov reflection method with respect to the integral leaves of
$\mathcal{D}$. That this makes sense we could see noticing that the
flux of $Y$ is, at fixed $s$, an ambient isometry.

We remark that the integrability condition on $\mathcal{D}$ imposes
that the form $\omega =\ae Y, \cdot \ad$ satisfies $\dd \omega=0$ on
$\mathcal{D}\times \mathcal{D}$. This implies that $\ae D_v Y,
w\ad=0$ for all vector fields $v,w$ on $\mathcal{D}$. So, the
integral leaves for $\mathcal{D}$ are totally geodesic on
$\mathbb{N}$.

Next, we use the height estimates on Theorem 11 to show the
existence of CMC Killing graphs for $\kappa\ge0$. We observe that if
$\Sigma$ is a Killing graph, then each flow line through $\Sigma$
meets $\Omega$. Since the norm of $Y$ is constant along the flow
lines, we have $\inf_\Sigma |Y|= \inf_\Omega |Y|$ and so on. Thus,
the height estimates on Theorem 11 for the particular case of graphs
depend only on data of the domain and the mean curvature.

Killing graphs corresponding to the vertical field $Y=\partial_t$
were previously studied in \cite{HLR} (see also \cite{AD}). We then
restrict ourselves to consider the horizontal field $Y=a\wedge
\partial_t\wedge p$.
Let $\Omega$ be a domain on $\Pi\cap \mathbb{N}$. We may write
$\Pi=\sigma\times\Rr$, where $\sigma$ is a horizontal geodesic
parametrized by $\rho$. Let $s=u(\rho,t)$ be a function defined on
$\Omega$, which specifies a point on the flow line of $Y$ starting
at the point with coordinates $(\rho,t)$ on $\Pi$. Let $\Sigma$ be
the (Killing) graph of $u$. We fix boundary data $u=0$ on
$\partial\Omega$. The tangent vectors to $\Sigma$ are
\[
X_\rho = \partial_\rho +\frac{\partial u}{\partial\rho}\,
\partial_s,\quad X_t =\partial_t + \frac{\partial u}{\partial t}\, \partial_s
\]
Recall that $\partial_s =Y$. The unit normal  vector field  to
$\Sigma$ so that $\ae Y,n\ad\ge 0$ is
\[
n = \frac{1}{W}\big(-\frac{\partial u}{\partial\rho}\,\partial_\rho
- \frac{\partial u}{\partial t}\,\partial_t + f\,\partial_s \big)
\]
where $W^2 =f+(\partial_\rho u)^2 +(\partial_t u)^2 = f+|\nabla
u|^2$. Thus
\[
\ae Y, n\ad = \frac{1}{W}.
\]
So, a lower estimate  for $\ae Y,n\ad$ gives an upper estimate for
$W$ and therefore for $|\nabla u|$. However, $\ae Y,n\ad$ is
super-harmonic for $\kappa\ge0$. Then
\[
\min_\Omega \ae Y,n\ad = \min_{\partial\Omega} \ae Y, n\ad.
\]
Then $|u|_{1,\Omega}\le C\, |u|_{1,\partial\Omega}$ for some
constant $C$ independent of $u$. We must then estimate $|\nabla u|$
on the boundary. As we said before, this means to get a lower
estimate for $\ae Y,n\ad$ on the boundary. Since $s=0$ on
$\partial\Omega$ and  $\phi$ is also super-harmonic, then
\[
\min_\Omega \ae Y,n\ad = \min_{\partial\Omega} \ae Y, n\ad
=\min_{\partial\Omega} \phi =\min_{\Omega} \phi=\phi(p_0)
\]
for a given $p_0$ on $\partial\Omega$. Thus, if $\eta$ is the unit
interior co-normal then at $p_0$
\begin{eqnarray*}
& & 0 \le \ae \nabla \phi, \eta\ad = H c \ae \nabla s, \eta\ad +\eta
\ae Y,n\ad = H c \,f \ae Y,\eta\ad +\ae D_\eta Y, n\ad +\ae Y,
D_\eta
n\ad\\
& & = Hc f\, \ae Y,\eta\ad + \ae a\wedge
\partial_t \wedge \eta, n\ad + \ae D_\eta n, \eta\ad \ae Y,\eta\ad\\
& & =Hc f\, \ae Y,\eta\ad + \ae \tau, n\ad - \ae A(\eta),\eta\ad\ae
Y,\eta\ad = \big(Hcf - \ae A\eta, \eta\ad\big)\, \ae Y,\eta\ad,
\end{eqnarray*}
where $\tau$ is the unit positively oriented vector tangent to
$\partial\Omega$. However by choice of $c$ we have $cf \le 1$ and
$cf\ge \gamma=:\inf_\Omega |Y|^2/\sup_\Omega |Y|^2$. Then $0\le
\big(H\gamma - \ae A\eta, \eta\ad\big)\, \ae Y,\eta\ad$. Since $\ae
Y,\eta\ad \ge 0$ and $2H=\ae A\eta, \eta \ad + \ae A\tau, \tau\ad$
then $\ae A\tau,\tau\ad \ge H(2-\gamma)$. However, $n=\ae n,
\bar\eta\ad \bar\eta +\ae n,fY\ad Y$ where $\bar\eta$ is the unit
interior normal to $\partial\Omega$ on $\Pi$. Then
\begin{eqnarray*}
& & \ae A\tau, \tau \ad = \ae D_\tau n, \tau\ad =\ae n, \bar\eta\ad
\ae D_\tau \bar\eta, \tau \ad +\ae n, fY \ad \ae D_{\tau}Y, \tau \ad
= \ae n, \bar\eta \ad \, \kappa_g,
\end{eqnarray*}
where $\kappa_g$ is the geodesic curvature on $\partial\Omega$,
which we suppose to be strictly positive. Since $H\le 0$, $\ae n,
\bar\eta \ad \le 0$ and $\kappa_g> 0$ then at $p_0$ we have
\[
H^2 \tilde\gamma^2\ge \kappa_g^2 \ae n, \bar\eta\ad^2 = \kappa_g^2
\big(1-f\,\ae n, Y\ad^2\big)
\]
where $\tilde\gamma=2-\gamma$. Denoting $\inf_\Omega|Y|=\widetilde
c$ we then obtain
\[
|\ae Y, n \ad| \ge \widetilde c \, \sqrt{\frac{\kappa_g^2
-H^2\tilde\gamma^2}{\kappa_g^2}}
\]
We then verified that there exists a non-trivial gradient estimate
for $u$ if we suppose $|H|< \kappa_g/\tilde\gamma$. By Schauder
theory on quasi-linear elliptic equations, we conclude that there
exists a CMC Killing graph on $\mathbb{N}$ with boundary on $\Pi$.

\vspace{0.3cm} \noindent {\bf Theorem 12. } {\it Let $\Pi$ be a
vertical plane on the Riemannian product $\Mm\times\Rr$, $\kappa\ge
0$, determined by an unit vector $a$ in $\mathbb{E}^4$. Let $\Omega$
be a domain on $\Pi$ which does not contain points of the axis
$\{\pm a\}\times\Rr$. If $|H|< \kappa_g/\tilde\gamma$, where
$\kappa_g$ is the geodesic curvature of $\partial\Omega$ in $\Pi$
and $\tilde\gamma=2-\sup_\Omega f/\inf_\Omega f$, then there exists
a surface (a Killing graph) with constant mean curvature $H$ and
boundary $\partial\Omega$.}

\vspace{0.3cm}

Notice that $Y=\partial_\theta$, where $\theta$ is the polar
coordinate centered at $r\, a$ as defined on Section 1. Thus
$|Y|=\sn_\kappa (\rho)$. Thus, a simple application to the flux
formula gives us the following area estimate:
\[
|H|\le  \frac{\max_{\partial\Omega}
\textrm{sn}_\kappa(\rho)}{\min_\Omega \sn_\kappa (\rho) }\,
\frac{|\partial \Omega|}{2|\Omega|}.
\]
This estimate also holds for surfaces in $\Mm\times\Rr$ satisfying
the condition that its boundary bounds a domain on a vertical plane
which does not contain singularities of $Y$.

Next, we fix $\varepsilon=-1$ and $\kappa\le 0$ . Let $\Sigma$ be a
CMC surface whose boundary is a geodesic circle on some horizontal
plane $\Pp$. Thus considering the Killing field $\partial_t$, the
function $\phi$ we defined above becomes $\phi = H t -\nu$, where
$\nu =\ae n, \partial_t\ad $. Then we have as before
\[
\Delta \phi = \big(2H^2 -|A|^2 +\kappa (1-\nu^2)\big)\, \nu
\]
We recall that
\[
|A|^2 = 2H^2 +2|\psi^1|^2, \quad 4|\psi^2|^2 =\kappa^2 (1-\nu^2)^2
\]
and since $\kappa\le 0$ and $1-\nu^2\le 0$ we have $2|\psi^2|=
\kappa  (1-\nu^2)$.  Replacing this above and assuming that
$|\psi^1|^2 -|\psi^2|\ge 0$ we have
\[
\Delta \phi =-2\big(|\psi^1|^2 -|\psi^2|\big)\, \nu\ge 0,
\]
since that we choose $n$ pointing upwards (which implies that $H\le
0$). By Stokes's theorem
\begin{eqnarray*}
& & -2 \int_\Sigma \big(|\psi^1|^2 -|\psi^2|\big)\, \nu\, \dd A =
\int_{\partial \Sigma}\ae \nabla \phi, \eta\ad \, \dd \sigma
\end{eqnarray*}
where $\eta$ is the outward unit co-normal to $\Sigma$ along
$\partial\Sigma$. However
\begin{eqnarray*}
& & \ae\nabla \phi,\eta\ad = H \ae\nabla t, \eta\ad -\ae\nabla \nu,
\eta\ad = -H\ae\partial_t, \eta\ad +\ae \partial_t, A\eta\ad
\end{eqnarray*}
Therefore, $\ae\nabla \phi,\eta\ad = \big(\ae A\eta,
\eta\ad-H\big)\,\ae \eta,
\partial_t\ad$. However, $\ae A\eta, \eta\ad=2H-\ae A\tau,\tau\ad$, where $\tau$ is
the unit tangent vector to $\partial\Sigma$. Let $\bar \eta$ be the
outwards unit normal to $\partial\Sigma$ with respect to $\Pp$.
Since $n= \ae n, \bar\eta\ad \bar\eta -\ae n,
\partial_t\ad \, \partial_t$ and since $\tau$ is orthogonal to both $\partial_t$
and $\bar\eta$ it follows that
\[
-\ae A \tau, \tau\ad = \ae D_{\tau}n, \tau\ad =\ae n, \bar\eta\ad
\ae D_{\tau}\bar\eta, \tau\ad = -\kappa_g \ae n,\bar\eta\ad
=\kappa_g \ae
\partial_t,\eta\ad
\]
Thus we conclude that
$\ae \nabla \phi, \eta \ad =\big(H +\kappa_g \ae \eta, \partial_t
\ad\big) \, \ae \eta, \partial_t\ad$. So by flux formula
\begin{eqnarray*}
& & \int_{\partial\Sigma}\, \ae \nabla \phi, \eta\ad \, \dd\sigma =H
\int_{\partial\Sigma} \ae \eta, \partial_t \ad \dd \sigma
+\int_{\partial\Sigma} \kappa_g \ae \eta, \partial_t \ad^2 \dd
\sigma =
2H^2 |\Omega| +\\
& &  \int_{\partial\Sigma} \kappa_g \ae \eta,
\partial_t \ad^2 \dd \sigma
\end{eqnarray*}
Gathering the expressions, we have
\[
-2\int_\Sigma \big(|\psi^1|^2 -|\psi^2|\big)\, \nu\, \dd A =2H^2
|\Omega| + \int_{\partial\Sigma} \kappa_g \ae \eta,
\partial_t \ad^2 \dd \sigma.
\]
Now again by flux formula
\[
\Big(\int_{\partial\Sigma}\ae \eta, \partial_t\ad \, \dd
\sigma\Big)^2 = 4H^2 |\Omega|^2
\]
But by Cauchy-Schwarz on $L^2$ functions we have
\[
\Big(\int_{\partial\Sigma}\ae \eta, \partial_t\ad \, \dd
\sigma\Big)^2\le |\partial\Sigma|\, \int_{\partial\Sigma}\ae \eta,
\partial_t\ad ^2 \, \dd \sigma
\]
So
\[
\frac{4H^2 |\Omega|^2}{|\partial \Sigma|}\le
\int_{\partial\Sigma}\ae \eta,
\partial_t\ad ^2 \, \dd \sigma
\]
Thus
\begin{eqnarray}
\label{angle} & & -2\int_\Sigma \big(|\psi^1|^2 -|\psi^2|\big)\,
\nu\, \dd A \le 2H^2
\frac{|\Omega|}{\partial\Omega}\big(|\partial\Omega|+2|\Omega|\,
\kappa_g\big)
\end{eqnarray}
with equality if and only of $\ae \eta, \partial_t\ad$ is constant
along $\partial\Sigma$. Now,  the geodesic curvature of
$\partial\Sigma$ calculated with respect to $\bar\eta$ is
$\kappa_g=-\ct_\kappa(\rho)$.  Thus
\begin{eqnarray*}
& & |\partial\Omega|+2|\Omega|\, \kappa_g  =\frac{2\pi}\kappa\,
\sn_\kappa(\rho)\, \big(\cs_\kappa(\rho)-1\big)^2\le 0
\end{eqnarray*}
since $\kappa\le 0$. So, occurs  equality on (\ref{angle}). Then,
the angle between $\Sigma$ and the horizontal plane is constant
along $\partial\Sigma$. So, $\Sigma$ is a stationary  surface for
the energy defined on Section 5. Thus, by Theorem 9,  $Q=0$ and the
surface is a hyperbolic cap.

\vspace{0.3cm} \noindent {\bf Theorem 13.} {\it Fix $\varepsilon=-1$
and $\kappa\le 0$. Let $\Sigma$ be a immersed CMC surface whose
boundary is a geodesic circle on some horizontal plane $\Pp$. If we
suppose that $|\psi^1|^2 - |\psi^2|\ge 0$, then $Q=0$ and the
surface is part of a hyperbolic cap or a planar disc.}

\vspace{0.3cm}

This theorem is a partial answer to a Lorentzian formulation of the
well-known {\it spherical cap conjecture} which was positively
proved on \cite{AP2} for the case $\kappa=0$.



\vspace{2cm}
\noindent Marcos Petr\'ucio Cavalcante\\
IMPA\\
Estrada Dona Castorina, 110\\
Rio de Janeiro, Brazil\\
22460-320\\
petrucio@impa.br\\

\vspace{0.5cm}

\noindent Jorge Herbert S. de Lira\\
Departamento de Matem\'atica - UFC\\
Campus do Pici, Bloco 914\\
Fortaleza, Cear\'a, Brazil\\
60455-760\\
jherbert@mat.ufc.br

\end{document}